# ON SHINTANI ZETA FUNCTIONS FOR GL(2)

Akihiko Yukie[1]

Oklahoma State University

ABSTRACT. In this paper we consider an analogue of the zeta function for not necessarily prehomogeneous representations of GL(2) and compute some of the poles.

**Table of contents.**

**Introduction**

The notion of prehomogeneous vector spaces was introduced by Sato and Shintani in [sash]. The main object of the theory is the zeta function, which is the counting function for generic rational orbits. The poles and the residues of the zeta functions are closely related to distributions of arithmetic objects. When the group is GL(2), Shintani determined the poles and the residues of the zeta functions for the spaces of binary forms of degree 2 and 3. His results, combined with the local theory of Datskovsky–Wright [dawra], [dawrb] and Datskovsky [dats] produced the zeta function theoretic proofs of Davenport–Heilbronn theorem on the density of cubic fields and Goldfeld–Hoffstein theorem on the density of class numbers of quadratic fields.

One basic reason why the zeta function theory yields such density theorems is that the orbit spaces of prehomogeneous vector spaces parametrize interesting arithmetic objects. However, as far as parametrizing arithmetic objects is concerned, we do not have to restrict ourselves to prehomogeneous vector spaces. For example, the orbit space of binary forms of degree 4 parametrizes rational isomorphism classes of elliptic curves defined over a number field.

The definition of the zeta function can easily be generalized to general representations of reductive groups using the notion of stability in geometric invariant theory. In this paper, we consider such a zeta function for the space of binary forms

---

[1]1991 Mathematics Subject Classification 11M41 The author is partially supported by NSF

of degree $d$. Our main theorem is Theorem (6.12). As the reader will see, our result does not determine all the poles of the zeta function. Also at present, we do not know how to deduce density theorems for orbits from the knowledge of the poles and the residues of the global zeta function.

## §1 Basic definitions

Throughout this paper, $k$ is a number field.

In this section, we consider the space of binary forms of degree $d \geq 4$, and define the zeta function and related notions. The cases $d = 2, 3$ were considered by Shintani, Wright, and the author (see [shintania], [shintanib], [wright], [yukiea]).

We basically follow the notations of [yukiec], but we recall the most basic ones. If $f, g$ are functions on a set $X$, $f \ll g$ means that there exists a constant $C$ such that $|f(x)| \leq Cg(x)$ for all $x \in X$. We also use the classical notation $x \ll y$ when $y$ is a much larger number than $x$. We hope the meaning of this notation will be clear from the context. The ring of adeles (resp. the group of ideles) over $k$ is denoted by $\mathbb{A}$ (resp. $\mathbb{A}^\times$). For a vector space $V$ over $k$, $V_\mathbb{A}$ is the adelization. We define $\mathbb{R}_+ = \{x \in \mathbb{R} \mid x > 0\}$. For $\lambda \in \mathbb{R}_+$, $\underline{\lambda}$ is the idele whose component at any infinite place is $\lambda^{\frac{1}{[k:\mathbb{Q}]}}$ and whose component at any finite place is 1. Let $|x|$ be the adelic absolute value of $x \in \mathbb{A}$. For any place $v$ of $k$, let $k_v$ be the completion of $k$ with respect to $v$, and $|x|_v$ the norm on $k_v$. Then $|\underline{\lambda}| = \lambda$. Let $\mathbb{A}^1 = \{t \in \mathbb{A}^\times \mid |t| = 1\}$. We normalize the measures $dt, d^\times t$ on $\mathbb{A}, \mathbb{A}^1$ so that the volumes of $\mathbb{A}/k$, $\mathbb{A}^1/k^\times$ are 1. Let $d\mu$ be the usual Haar measure on $\mathbb{R}$. We use the measure $d^\times \mu = \mu^{-1} d\mu$ as the measure on $\mathbb{R}_+$. We use the measure $d^\times \mu d^\times t$ for $(\mu, t) \in \mathbb{R}_+ \times \mathbb{A}^1 \cong \mathbb{A}^\times$. We fix a non-trivial character $\langle \ \rangle$ on $\mathbb{A}/k$.

Let $V = \text{Sym}^d k^2$ be the vector space of degree $d$ forms in two variables $v = (v_1, v_2)$. When $d$ is even, we put $d = 2d_1$. The vector space $V$ can be identified with $k^{d+1}$ by the following map

$$f(v) = f_x(v_1, v_2) = x_0 v_1^d + \cdots + x_d v_2^d \to (x_0, \cdots, x_d).$$

The group $\text{GL}(2)$ acts on $V$ from the left by the formula

$$(g \cdot f)(v) = f(vg).$$

Let $G_1 = \text{GL}(1), G_2 = \text{GL}(2)$, and $G = G_1 \times G_2$. We consider the ordinary multiplication by elements of $G_1$ to $V$. This defines an action of $G$ on $V$.

Let $T_2 \subset G_2$ be the subgroup of diagonal matrices, and $N_2 \subset G_2$ the subgroup of lower triangular matrices with diagonal entries 1. We define $T = G_1 \times T_2$ and $N = \{1\} \times N_2$. Then $B = TN$ is a Borel subgroup of $G$. For $i = 0, \cdots, d$, let $\gamma_i$ be the weight of $x_i$, i.e. the $i$-th coordinate of $t(x_0, \cdots, x_d)$ is $\gamma_i(t) x_i$ for $t \in T_\mathbb{A}$.

We write elements of $T_2, N_2$ in the following manner:

$$a(t) = a(t_1, t_2) = \begin{pmatrix} t_1 & 0 \\ 0 & t_2 \end{pmatrix}, \quad n(u) = \begin{pmatrix} 1 & 0 \\ u & 1 \end{pmatrix}.$$

We identify $n(u), a(t)$ with $(1, n(u)), (1, a(t)) \in B_\mathbb{A}$ and use the same notation. For $t_1, t_2 \in \mathbb{A}^\times$ and $z = (z_1, z_2) \in \mathbb{C}^2$, we define $a(t_1, t_2)^z = |t_1|^{z_1}|t_2|^{z_2}$. Let

$$G_{2\mathbb{A}}^0 = \{g_2 \in G_{2\mathbb{A}} \mid |\det g_2| = 1\},\ G_\mathbb{A}^0 = \mathbb{A}^1 \times G_{2\mathbb{A}}^0,$$
$$B_{2\mathbb{A}}^0 = G_{2\mathbb{A}}^0 \cap B_{2\mathbb{A}},\ B_\mathbb{A}^0 = \mathbb{A}^1 \times B_{2\mathbb{A}}^0,$$
$$\widehat{B}_{2\mathbb{A}}^0 = \{a(t_1, t_2)n(u) \mid t_1, t_2 \in \mathbb{A}^1, u \in \mathbb{A}\},\ \widehat{B}_\mathbb{A}^0 = \mathbb{A}^1 \times \widehat{B}_{2\mathbb{A}}^0,$$
$$T_\mathbb{A}^0 = T_\mathbb{A} \cap B_\mathbb{A}^0,\ \widehat{T}_\mathbb{A}^0 = T_\mathbb{A} \cap \widehat{B}_\mathbb{A}^0.$$

We define $\widetilde{G}_\mathbb{A} = \mathbb{A}^\times \times G_{2\mathbb{A}}^0 \cong \mathbb{R}_+ \times G_\mathbb{A}^0$ and $\widetilde{B}_\mathbb{A} = \mathbb{R}_+ \times B_\mathbb{A}^0$. The group $\widetilde{G}_\mathbb{A}$ acts on $V_\mathbb{A}$ also. Throughout this paper, we express elements $\widetilde{g} \in \widetilde{G}_\mathbb{A}$, $g^0 \in G_\mathbb{A}^0$ as $\widetilde{g} = (\lambda, t, g_2)$, $g^0 = (t, g_2)$, where $\lambda \in \mathbb{R}_+$, $t \in \mathbb{A}^1$, and $g_2 \in G_{2\mathbb{A}}^0$.

Let $K$ be the usual maximal compact subgroup of $G_\mathbb{A}$ i.e.

$$K \cap G_{2\mathbb{A}} = \prod_{v \in \mathfrak{M}_\mathbb{R}} \mathrm{O}(2) \times \prod_{v \in \mathfrak{M}_\mathbb{C}} \mathrm{U}(2) \times \prod_{v \in \mathfrak{M}_f} \mathrm{GL}(2, o_v)$$

($\mathbb{A}^\times$ has a unique maximal compact subgroup). Let $dk$ be the measure on $K$ such that $\int_K dk = 1$. Throughout this paper, we express elements of $T_\mathbb{A}^0$, $\widehat{T}_\mathbb{A}^0$, $B_\mathbb{A}^0$, $\widehat{B}_\mathbb{A}^0$, and $\widetilde{B}_\mathbb{A}$ as

$$t^0 = (t, a(\underline{\mu}t_1, \underline{\mu}^{-1}t_2)),\ \widehat{t}^{\,0} = (t, a(t_1, t_2)),$$
$$b^0 = t^0 n(u) = n(u_0)t^0,\ \widehat{b}^{\,0} = \widehat{t}^{\,0} n(u) = n(\widehat{u}_0)\widehat{t}^{\,0},\ \text{and}\ \widetilde{b} = (\lambda, b^0),$$

where $t, t_1, t_2 \in \mathbb{A}^1, \mu \in \mathbb{R}_+, u \in \mathbb{A}$, and $u_0 = \underline{\mu}^{-2}t_1^{-1}t_2 u, \widehat{u}_0 = t_1^{-1}t_2 u$ respectively. The group $G_\mathbb{A}^0$ has the Iwasawa decomposition $G_\mathbb{A}^0 = KT_\mathbb{A}^0 N_\mathbb{A}$. Let

$$g^0 = k(g^0)(t(g^0), a(t(g^0)))n(u(g^0)) = k(g^0)(t(g^0), a(t_1(g^0), t_2(g^0)))n(u(g^0))$$

be the Iwasawa decomposition of the element $g^0 \in G_\mathbb{A}^0$.

We use the measures

$$d^\times \widehat{t}^{\,0} = d^\times t d^\times t_1 d^\times t_2,\ d^\times t^0 = d^\times \mu d^\times \widehat{t}^{\,0},$$
$$d\widehat{b}^{\,0} = d^\times \widehat{t}^{\,0} du,\ db^0 = \mu^{-2} d^\times t^0 du = d^\times t^0 du_0,\ \text{and}\ d\widetilde{b} = d^\times \lambda db^0$$

on $\widehat{T}_\mathbb{A}^0$, $T_\mathbb{A}^0$, $\widehat{B}_\mathbb{A}^0$, $B_\mathbb{A}^0$, and $\widetilde{B}_\mathbb{A}$ respectively. We use the measure $dg^0 = dk db^0$ as the Haar measure on $G_\mathbb{A}^0$. Let $dg_2$ be the Haar measure on $G_{2\mathbb{A}}^0$ such that $dg^0 = d^\times t dg_2$, where $d^\times t$ is the measure on $\mathbb{A}^1$ we defined earlier. We use $d\widetilde{g} = d^\times \lambda dg^0$ as the measure on $\widetilde{G}_\mathbb{A}$.

Let $\zeta_k(s)$ be the Dedekind zeta function. As in [weilc], we define

$$Z_k(s) = |\Delta_k|^{\frac{s}{2}} \left(\pi^{-\frac{s}{2}} \Gamma\left(\frac{s}{2}\right)\right)^{r_1} \left((2\pi)^{-s} \Gamma(s)\right)^{r_2} \zeta_k(s),$$

where $\Delta_k$ is the discriminant of $k$. Let $\phi(s) = \frac{Z_k(s)}{Z_k(s+1)}$. We define

$$\mathfrak{R}_k = \operatorname*{Res}_{s=1} Z_k(s),\ \varrho = \operatorname*{Res}_{s=1} \phi(s) = \frac{\mathfrak{R}_k}{Z_k(2)},\ \text{and}\ \mathfrak{V}_2 = \varrho^{-1}.$$

The volume of $G_{\mathbb{A}}^0/G_k$ with respect to the measure $dg^0$ is $\mathfrak{V}_2$. For a character $\omega$ of $\mathbb{A}^\times/k^\times$, we define $\delta(\omega) = 1$ if $\omega$ is trivial, and $\delta(\omega) = 0$ otherwise.

We defined the notion of $k$-stable points in Definition (3.1.3) in [yukiec] only for prehomogeneous vector spaces. However, the definition applies to general representations without any modification. For example, in our case, a point $x \in V_k$ is semi-stable (resp. $k$-stable) if it has no factor of multiplicity $> \left[\frac{n}{2}\right]$ (resp. no $k$-rational factor of multiplicity $\geq \left[\frac{n}{2}\right]$). Note that if $x$ has a factor of multiplicity $> \left[\frac{n}{2}\right]$, that factor is rational over $k$. We denote the set of semi-stable points (resp. $k$-stable points) by $V_k^{\mathrm{ss}}$ (resp. $V_k^{\mathrm{s}}$). It is clear from this definition that $V_k^{\mathrm{ss}} = V_k^{\mathrm{s}}$ if $n$ is odd. When $d$ is even, we define $V_{\mathrm{st}k}^{\mathrm{ss}} = V_k^{\mathrm{ss}} \setminus V_k^{\mathrm{s}}$ and call it the set of strictly semi-stable points.

For any vector space $W$ over $k$, let $\mathscr{S}(W_{\mathbb{A}})$, $\mathscr{S}(W_{\mathbb{A}_f})$, and $\mathscr{S}(W_{\mathbb{A}_\infty})$ be the space of Schwartz–Bruhat functions. For $\Phi \in \mathscr{S}(V_{\mathbb{A}})$, we define

$$(1.1) \qquad \Theta_V^{\mathrm{s}}(\Phi, \widetilde{g}) = \sum_{x \in V_k^{\mathrm{s}}} \Phi(\widetilde{g}x), \; \Theta_{V,\mathrm{st}}(\Phi, \widetilde{g}) = \sum_{x \in V_{\mathrm{st}k}^{\mathrm{ss}}} \Phi(\widetilde{g}x).$$

Let $\omega = (\omega_1, \omega_2)$ be a character of $(\mathbb{A}^\times/k^\times)^2$. For $\widetilde{g} = (\lambda, t, g_2) \in \widetilde{G}_{\mathbb{A}}$, we define $\omega(\widetilde{g}) = \omega_1(t)\omega_2(\det g_2)$. This defines a character of $\widetilde{G}_{\mathbb{A}}/G_k$.

**Definition (1.2)** *Let $\omega$ be as above. For $\Phi \in \mathscr{S}(V_{\mathbb{A}})$ and a complex variable $s$, we define*

$$(1) \qquad Z(\Phi, \omega, s) = \int_{\widetilde{G}_{\mathbb{A}}/G_k} \lambda^s \omega(\widetilde{g}) \Theta_V^{\mathrm{s}}(\Phi, \widetilde{g}) d\widetilde{g},$$

$$(2) \qquad Z_+(\Phi, \omega, s) = \int_{\substack{\widetilde{G}_{\mathbb{A}}/G_k \\ \lambda \geq 1}} \lambda^s \omega(\widetilde{g}) \Theta_V^{\mathrm{s}}(\Phi, \widetilde{g}) d\widetilde{g}.$$

By the proof of Proposition (3.1.4) in [yukiec], (1) converges absolutely for $\mathrm{Re}(s) \gg 0$, and (2) is an entire function.

Let $\binom{d}{i}$ be the binomial coefficient. We define a non-degenerate bilinear form $[\,,\,]'$ on $V_{\mathbb{A}}$ by the formula

$$[x, y]' = \sum_{i=0}^d \binom{d}{i}^{-1} x_i y_i$$

for $x = (x_0, \cdots, x_d)$ and $y = (y_0, \cdots, y_d)$. Easy computations show that $[g_2 x, y]' = [x, {}^t g_2 y]'$ for all $x, y$. Let $\nu = \begin{pmatrix} 0 & 1 \\ 1 & 0 \end{pmatrix}$. We define

$$[x, y] = [x, \nu y]', \; \widetilde{g}^\iota = (\lambda^{-1}, t^{-1}, \nu^{-1 t} g_2^{-1} \nu).$$

Then $[\widetilde{g}x, y] = [x, (\widetilde{g}^\iota)^{-1} y]$ for all $x, y \in V_{\mathbb{A}}$. For $\Phi \in \mathscr{S}(V_{\mathbb{A}})$, we define the Fourier transform $\widehat{\Phi}$ of $\Phi$ by the formula

$$\widehat{\Phi}(x) = \int_{V_{\mathbb{A}}} \Phi(y) \langle [x, y] \rangle dy.$$

It is easy to see that the Fourier transform of the function $\Phi(\widetilde{g}\cdot)$ is $\lambda^{-(d+1)}\widehat{\Phi}(\widetilde{g}^\iota\cdot)$.

For $\lambda \in \mathbb{R}_+$, we define $\Phi_\lambda(x) = \Phi(\underline{\lambda}x)$.

**Definition (1.3)** *For $\Phi, \omega$, and $s$ as above, we define*

(1) $\qquad I(\Phi, \omega) = \int_{G_\mathbb{A}^0/G_k} \omega(g^0) \left( \sum_{x \in V_k \setminus V_k^s} \widehat{\Phi}((g^0)^\iota x) - \sum_{x \in V_k \setminus V_k^s} \Phi(g^0 x) \right) dg^0,$

(2) $\qquad I(\Phi, \omega, s) = \int_0^1 \lambda^s I(\Phi_\lambda, \omega) d^\times \lambda.$

By the Poisson summation formula,

$$Z(\Phi, \omega, s) = Z_+(\Phi, \omega, s) + Z_+(\widehat{\Phi}, \omega^{-1}, d+1-s) + I(\Phi, \omega, s).$$

For $i > [\frac{d}{2}]$, we define

$$Y_{ik} = \{x = (\overbrace{0, \cdots, 0}^{i}, x_i, \cdots, x_d) \mid x_i, \cdots, x_d \in k\},$$
$$Y_{ik}^{ss} = \{x \in Y_{ik} \mid x_i \in k^\times\},$$
$$Z_{ik} = \{x = (\overbrace{0, \cdots, 0}^{i}, x_i, \overbrace{0, \cdots, 0}^{d-i-1}) \mid x_i \in k\},$$
$$Z_{ik}^{ss} = \{x \in Z_{ik} \mid x_i \in k^\times\}.$$

Let $S_{ik} = G_k Y_{ik}^{ss}$ for $i > [\frac{d}{2}]$. Then it is easy to see that $S_{ik} \cong G_k \times_{B_k} Y_{ik}^{ss}$.

Suppose $d$ is even. For $d_1 + 2 \leq i \leq d$, we define

$$Y_{st,ik} = \{(\overbrace{0, \cdots, 0}^{d_1}, x_{d_1}, \overbrace{0, \cdots, 0}^{i-d_1-1}, x_i, \cdots, x_d) \mid x_{d_1}, x_i, \cdots, x_d \in k\},$$
$$Y_{st,ik}^{ss} = \{x \in Y_{st,ik} \mid x_{d_1}, x_i \in k^\times\},$$
$$Z_{st,ik} = \{(\overbrace{0, \cdots, 0}^{d_1}, x_{d_1}, \overbrace{0, \cdots, 0}^{i-d_1-1}, x_i, \overbrace{0, \cdots, 0}^{d-i}) \mid x_i \in k\},$$
$$Z_{st,ik}^{ss} = \{x \in Z_{st,ik} \mid x_{d_1}, x_i \in k^\times\},$$
$$Z'_{0k} = \{(\overbrace{0, \cdots, 0}^{d_1}, x_{d_1}, \overbrace{0, \cdots, 0}^{d_1}) \mid x_{d_1} \in k\},$$
$$Z'^{ss}_{0k} = \{x \in Z'_{0k} \mid x_{d_1} \in k^\times\}.$$

Let $S_{st,ik} = G_k Y_{st,ik}^{ss}$ for $d_1 + 2 \leq i \leq d$. Let $H$ the subgroup of $G$ generated by $T$ and $\nu$. Then it is easy to see that $V_{stk}^{ss} \cong G_k \times_{H_k} Z'^{ss}_{0k}$ and $S_{st,ik} \cong G_k \times_{T_k} Y_{st,ik}^{ss}$.

The following lemma is easy to verify, and the proof is left to the reader.

**Lemma (1.4)**
(1) If $d$ is odd, $V_k \setminus \{0\} = V_k^s \coprod_{[\frac{d}{2}] < i \leq d} S_{ik}$.
(2) If $d = 2d_1$, $V_k \setminus \{0\} = V_k^s \coprod V_{stk}^{ss} \coprod_{d_1 < i \leq d} S_{ik} \coprod_{d_1 + 2 \leq i \leq d} S_{st,ik}$.

For $\Phi \in \mathscr{S}(V_\mathbb{A})$ and $\widetilde{g} \in \widetilde{G}$, we define

(1.5)
$$\Theta_{Y_i}(\Phi, \widetilde{g}) = \sum_{x \in Y_{ik}^{\mathrm{ss}}} \Phi(\widetilde{g}x),$$
$$\Theta_{S_i}(\Phi, \widetilde{g}) = \sum_{x \in S_{ik}} \Phi(\widetilde{g}x),$$
$$\Theta_{Z_0'}(\Phi, \widetilde{g}) = \sum_{x \in Z_{0k}'^{\mathrm{ss}}} \Phi(\widetilde{g}x).$$

We define $\Theta_{Y_{\mathrm{st},i}}(\Phi, \widetilde{g}) \Theta_{S_{\mathrm{st},i}}(\Phi, \widetilde{g})$ similarly.

For later purposes, we define an operator $M_\omega$ which is similar to the one in [wright].

**Definition (1.6)** *For $\Phi \in \mathscr{S}(V_\mathbb{A})$ and $\omega$ as above, we define*

$$M_\omega \Phi(x) = \int_K \omega(\det k) \Phi(kx) dk.$$

The operator $M_\omega \Phi(x)$ satisfies similar properties to those in Lemma 5.1 in [wright].

It is easy to see that $Z(\Phi, \omega, s) = Z(M_\omega \Phi, \omega, s)$. So for the rest of this paper, we assume that $\Phi = M_\omega \Phi$.

### §2 Estimates of theta series

In this section, we define and estimate various theta series. Let $C \subset \mathbb{A}^1$ be a compact subset which surjects to $\mathbb{A}^1/k^\times$. Throughout this section, we consider $\widetilde{b}$ such that $t, t_1, t_2 \in C$.

Let $i$ be a positive integer. For any $N > 0$ and $\lambda = (\lambda_1, \cdots, \lambda_i) \in \mathbb{R}_+^i$, we define

$$\mathrm{rd}_{i,N}(\lambda) = \inf(\lambda_1^{\pm N} \cdots \lambda_i^{\pm N}),$$

where we consider all the possible $\pm$.

Let $[\frac{d}{2}] \leq i < j \leq d$, and $x = (\overbrace{0, \cdots, 0}^{i}, x_i, \cdots, x_d)$. Since $n(u)$ is a lower triangular matrix, $n(u)x$ is of the form

$$(\overbrace{0, \cdots, 0}^{i}, x_i, x_{i+1} + p_{ii+1}, \cdots, x_d + p_{id}),$$

where $p_{ij} = p_{ij}(x, u) \in \mathbb{A}$ does not depend on $x_j, \cdots, x_d$. So

$$\widetilde{b}x = (\overbrace{0, \cdots, 0}^{i}, \underline{\lambda}\gamma_i(t^0)x_i, \underline{\lambda}\gamma_{i+1}(t^0)(x_{i+1} + p_{ii+1}), \cdots, \underline{\lambda}\gamma_d(t^0)(x_d + p_{id})).$$

We define

$$A_{ij}(\widetilde{b}, x) = (\underline{\lambda}\gamma_i(t^0)x_i, \underline{\lambda}\gamma_{i+1}(t^0)(x_{i+1} + p_{ii+1}), \cdots, \underline{\lambda}\gamma_{j-1}(t^0)(x_{j-1} + p_{ij-1})).$$

In other words,
$$\widetilde{b}x = (\overbrace{0,\cdots,0}^{i}, A_{ij}(\widetilde{b},x), \overbrace{\cdots}^{d-j+1}).$$
$(A_{ii+1}(\widetilde{b},x) = \underline{\lambda}\gamma_i(t^0)x_i.)$

Let
$$Y_{ijk} = \{(\overbrace{0,\cdots,0}^{i}, x_i, \cdots, x_{j-1}, \overbrace{0,\cdots,0}^{d-j+1}) \mid x_i, \cdots, x_{j-1} \in k\},$$
$$Y'_{ijk} = \{x \in Y_{ijk} \mid x_i \in k^{\times}\}.$$

**Definition (2.1)** *Let $\frac{d}{2} \leq i < j \leq d$. For $\Psi \in \mathscr{S}(\mathbb{A}^{j-i+1})$, we define*
$$\Theta_{ij}(\Psi, \widetilde{b}) = \sum_{\substack{x \in Y'_{ijk} \\ x_j \in k^{\times}}} \Psi(A_{ij}(\widetilde{b},x), \underline{\lambda}^{-1}\gamma_j(t^0)^{-1}x_j)\langle -x_j p_{ij}(x,u)\rangle.$$

Note that $\Theta_{ij}(\Psi, \widetilde{b})$ is a function on $\widetilde{B}_{\mathbb{A}}/B_k$.

**Proposition (2.2)** *For any $N \geq 1$,*
$$|\Theta_{ij}(\Psi,\widetilde{b})| \ll \begin{cases} \lambda^{-(j-i+N-2)}\mu^{(2i-d)N-(2j-d)} & \text{if } \lambda \leq 1, \mu \leq 1, \\ \mathrm{rd}_{2,N}(\lambda,\mu) & \text{otherwise.} \end{cases}$$

*Proof.* By Lemma (1.2.8) in [yukiec], for any $N_1, N_2 \geq 1$,
$$|\Theta_{ij}(\Psi,\widetilde{b})| \ll (\lambda\mu^{d-2i})^{-N_1}(\lambda\mu^{d-2j})^{N_2}\prod_{l=i+1}^{j-1}\sup(1, \lambda^{-1}\mu^{2l-d}).$$

If $\lambda \geq 1, \mu \geq 1$, choose $N_1 \gg N_2 \gg 0$ and $(2i-d)N_1 \ll (2j-d)N_2$. This is possible because $2j-d > 2i-d > 0$. If $\lambda \geq 1, \mu \leq 1$, choose $N_1 \gg 0$. If $\lambda \leq 1, \mu \geq 1$, choose $N_2 \gg 0$. If $\lambda, \mu \leq 1$, choose $N_2 = 1$. Then, since $\lambda^{-1}\mu^{2l-d} \leq \lambda^{-1}$ for $l = i+1, \cdots, j-1$ and $\lambda^{-1} \geq 1$, we get our estimate. $\square$

Next, we define another type of theta series when $d = 2d_1$. Let $d_1 + 2 \leq i \leq d$, and
$$x = (\overbrace{0,\cdots,0}^{d_1}, x_{d_1}, \overbrace{0,\cdots,0}^{i-d_1-1}, x_i, \cdots, x_d).$$

As before, $n(u)x$ is of the form
$$(\overbrace{0,\cdots,0}^{d_1}, x_{d_1}, c_{d_1+1}x_{d_1}u, \cdots, c_{i-1}x_{d_1}u^{i-d_1-1}, x_i + q_{ii}, \cdots, x_d + q_{id}),$$

where $c_l = \binom{d_1}{l-d_1}$, $q_{ii} = c_i x_{d_1} u^{i-d_1}$, and $q_{ir} = q_{ir}(x,u) \in \mathbb{A}$ does not depend on $x_r, \cdots, x_d$. Note that $c_{d_1+1} = d_1$.

Let
$$A'_{\mathrm{st},i}(\widetilde{b},x) = (\underline{\lambda}\gamma_{d_1}(t^0)x_{d_1}, \underline{\lambda}\gamma_{d_1+1}(t^0)c_{d_1+1}x_{d_1}u, \cdots, \underline{\lambda}\gamma_{i-1}(t^0)c_{i-1}x_{d_1}u^{i-d_1-1}).$$

Then
$$\widetilde{b}x = (\overbrace{0,\cdots,0}^{d_1}, A'_{\mathrm{st},i}(\widetilde{b},x), \overbrace{\cdots}^{d-i+1}).$$

For convenience, we define
$$A'_{\mathrm{st},d+1}(\widetilde{b},x) = (\underline{\lambda}\gamma_{d_1}(t^0)x_{d_1}, \underline{\lambda}\gamma_{d_1+1}(t^0)c_{d_1+1}x_{d_1}u, \cdots, \underline{\lambda}\gamma_d(t^0)c_d x_{d_1} u^{d_1}).$$

Note that $A'_{\mathrm{st},i}(\widetilde{b},x)$ depends only on $\widetilde{b}, x_{d_1}$ for all $i$.

For $d_1 + 2 \leq i < j \leq d$, we define
$$A_{\mathrm{st},ij}(\widetilde{b},x) = (A'_{\mathrm{st},i}(\widetilde{b},x), \underline{\lambda}\gamma_i(t^0)(x_i + q_{ii}), \cdots, \underline{\lambda}\gamma_{j-1}(t^0)(x_{j-1} + q_{ij-1})).$$

Let
$$Y_{\mathrm{st},ijk} = \{(\overbrace{0,\cdots,0}^{d_1}, x_{d_1}, \overbrace{0,\cdots,0}^{i-d_1-1}, x_i, \cdots, x_{j-1}, \overbrace{0,\cdots,0}^{d-j+1}) \mid x_{d_1}, x_i, \cdots, x_{j-1} \in k\},$$
$$Y'_{\mathrm{st},ijk} = \{x \in Y_{\mathrm{st},ijk} \mid x_{d_1}, x_i \in k^\times\}.$$

**Definition (2.3)** *For $\Psi \in \mathscr{S}(\mathbb{A}^{j-d_1+1})$, we define*
$$\Theta_{\mathrm{st},ij}(\Psi,\widetilde{b}) = \sum_{\substack{x \in Y'_{\mathrm{st},ijk} \\ x_j \in k^\times}} \Psi(A_{\mathrm{st},ij}(\widetilde{b},x), \underline{\lambda}^{-1}\gamma_j(t^0)^{-1}x_j)\langle -x_j q_{ij}(x,u)\rangle.$$

Note that $\Theta_{\mathrm{st},ij}(\Psi,\widetilde{b})$ is a function on $\widetilde{B}_\mathbb{A}/T_k$.

**Proposition (2.4)** *There exists a Schwartz–Bruhat function $0 \leq \Psi' \in \mathscr{S}(\mathbb{A}^2)$ such that for any $N_1, N_2 \geq 1$, $\Theta_{\mathrm{st},ij}(\Psi,\widetilde{b})$ is bounded by a constant multiple of*

$$(\lambda^2 \mu^{d-2i})^{-N_1}(\lambda \mu^{d-2j})^{N_2} \prod_{l=i+1}^{j-1} \sup(1, \lambda^{-1}\mu^{2l-d}) \sum_{x_{d_1} \in k^\times} \Psi'(\underline{\lambda}x_{d_1}, \underline{\lambda\mu}^{-2}x_{d_1}u).$$

*Proof.* Let
$$\overline{A}_{\mathrm{st},ij}(\lambda,\mu,u,x) = (\underline{\lambda\mu}^{d-2i}(x_i + c_i x_{d_1} u^{i-d_1}), \underline{\lambda\mu}^{d-2i-2}x_{i+1}, \cdots, \underline{\lambda\mu}^{d-2j+2}x_{j-1}).$$

By Proposition (1.2.3) and the proof of Lemma (1.2.7) in [yukiec], there exists a Schwartz–Bruhat function $0 \leq \Psi_1 \in \mathscr{S}(\mathbb{A}^2)$, such that $\Theta_{\mathrm{st},ij}(\Psi,\widetilde{b})$ is bounded by a constant multiple of

$$\sum_{\substack{x \in Y'_{\mathrm{st},ijk} \\ x_j \in k^\times}} \Psi_1(A'_{\mathrm{st},i}(\widetilde{b},x), \overline{A}_{\mathrm{st},ij}(\lambda,\mu,u,x), \underline{\lambda}^{-1}\underline{\mu}^{2j-d}x_j).$$

There exist Schwartz–Bruhat functions
$$0 \leq \Psi_2 \in \mathscr{S}(\mathbb{A}^{i-d_1+1}),\ 0 \leq \Psi_3 \in \mathscr{S}(\mathbb{A}^{j-i-1}),\ \text{and}\ 0 \leq \Psi_4 \in \mathscr{S}(\mathbb{A})$$

such that
$$\Psi_1(x_{d_1}, \cdots, x_j) \ll \Psi_2(x_{d_1}, \cdots, x_i)\Psi_3(x_{i+1}, \cdots, x_{j-1})\Psi_4(x_j).$$

Let
$$h_{N_2}(\lambda, \mu) = (\lambda \mu^{d-2j})^{N_2} \prod_{l=i+1}^{j-1} \sup(1, \lambda^{-1}\mu^{2l-d}).$$

By Lemma (1.2.8) in [yukiec], for any $N_2 \geq 1$, $\Theta_{\mathrm{st},ij}(\Psi, \widetilde{b})$ is bounded by a constant multiple of
$$h_{N_2}(\lambda, \mu) \sum_{x_{d_1}, x_i \in k^\times} \Psi_2(A'_{\mathrm{st},i}(\widetilde{b}, x), \underline{\lambda}\underline{\mu}^{d-2i}(x_i + c_i x_{d_1} u^{i-d_1})).$$

Let $B_{\mathrm{st},i}(\lambda, \mu, u, x)$ be the following vector
$$(\underline{\lambda}x_{d_1}, \underline{\lambda}\underline{\mu}^{-2}d_1 x_{d_1} u, \underline{\lambda}\underline{\mu}^{d-2i+2} c_{i-1} x_{d_1} u^{i-d_1-1}, \underline{\lambda}\underline{\mu}^{d-2i}(x_i + c_i x_{d_1} u^{i-d_1})).$$

Note that $B_{\mathrm{st},i}(\lambda, \mu, u, x)$ depends only on $\lambda, \mu, u, x_{d_1}, x_i$. We denote the finite part and the infinite part of $B_{\mathrm{st},i}(\lambda, \mu, u, x)$ by
$$B_{\mathrm{st},i,f}(\lambda, \mu, u, x), \ B_{\mathrm{st},i,\infty}(\lambda, \mu, u, x)$$
respectively.

There exists a Schwartz–Bruhat function $0 \leq \Psi_5 \in \mathscr{S}(\mathbb{A}^4)$ such that
$$\Psi_2(A'_{\mathrm{st},i}(\widetilde{b}, x), \underline{\lambda}\underline{\mu}^{d-2i}(x_i + c_i x_{d_1} u^{i-d_1})) \ll \Psi_5(B_{\mathrm{st},i}(\lambda, \mu, u, x)).$$

So we consider
$$\sum_{x_{d_1}, x_i \in k^\times} \Psi_5(B_{\mathrm{st},i}(\lambda, \mu, u, x)).$$

There exist Schwartz–Bruhat functions $0 \leq \Psi_6 \in \mathscr{S}(\mathbb{A}^4)$ and $0 \leq \Psi_7 \in \mathscr{S}(\mathbb{A}^2)$ such that
$$\Psi_5(x_{d_1}, x_{d_1+1}, x_{i-1}, x_i) \ll \Psi_6(x_{d_1}, x_{d_1+1}, x_{i-1}, x_i)\Psi_7(x_{d_1}, d_1^{-1}x_{d_1+1}).$$

So we only have to prove that for any $N_1 \geq 1$,

(2.5) $$\Theta(\Psi_6, \underline{\lambda}, \underline{\mu}, u, x_{d_1}) \stackrel{\text{def}}{=} \sum_{x_i \in k^\times} \Psi_6(B_{\mathrm{st},i}(\lambda, \mu, u, x))$$
$$\ll (\lambda^2 \mu^{d-2i})^{-N_1}.$$

We fix $x_{d_1} \in k^\times$. There exists an open compact subgroup $U_1 \subset \mathbb{A}_f$ such that
$$\Psi_6(x_{d_1}, x_{d_1+1}, x_{i-1}, x_i) = 0$$

unless $x_{lf} \in U_1$ for $l = d_1, d_1+1, i-1, i$ ($x_{lf}$ is the finite part of $x_l$). There exists an open compact subgroup $U_2 \subset \mathbb{A}_f$, such that $y, z \in U_1$ implies $yz, d_1^{-1} c_i c_{i-1}^{-1} yz \in U_2$. Let $L_1 = U_1 \cap k$, $L_2 = U_2 \cap k$. Then in the above sum, we only have to consider terms which satisfy

$$x_{d_1}, \ d_1 x_{d_1} u_f, \ c_{i-1} x_{d_1} u_f^{i-d_1-1}, \ x_i + c_i x_{d_1} u_f^{i-d_1} \in U_1.$$

($u_f$ is the finite part of $u$.)

But then

$$x_{d_1} x_i = x_{d_1}(x_i + c_i x_{d_1} u_f^{i-d_1}) - d_1^{-1} c_i c_{i-1}^{-1} d_1 x_{d_1} u_f c_{i-1} x_{d_1} u_f^{i-d_1-1}$$
$$\in (U_2 \setminus \{0\}) \cap k = L_2 \setminus \{0\}.$$

Therefore there exists a Schwartz–Bruhat function $0 \leq \Psi_{6\infty} \in \mathscr{S}(\mathbb{A}_\infty^4)$ such that

$$\Theta(\Psi_6, \underline{\lambda}, \underline{\mu}, u, x_{d_1}) \ll \sum_{x_{d_1} x_i \in L_2 \setminus \{0\}} \Psi_{6\infty}(B_{\mathrm{st},i,\infty}(\lambda, \mu, u, x)).$$

Let $\mathfrak{M}_\mathbb{R}, \mathfrak{M}_\mathbb{C}$ be the set of real places, imaginary places respectively. For $x \in \mathbb{A}_\infty$, we define

$$\|x\|_\infty = \sum_{v \in \mathfrak{M}_\mathbb{R}} |x|_v^{[k:\mathbb{Q}]} + \sum_{v \in \mathfrak{M}_\mathbb{C}} |x|_v^{\frac{[k:\mathbb{Q}]}{2}}.$$

Since $\Psi_{6\infty}$ is rapidly decreasing, for any $N_1 > 1$,

$$(1 + \|x_{d_1} x_i - d_1^{-1} c_i c_{i-1}^{-1} x_{d_1+1} x_{i-1}\|_\infty)^{N_1} \Psi_{6\infty}(x_{d_1}, x_{d_1+1}, x_{i-1}, x_i)$$

is bounded. So the above sum is bounded by a constant multiple of

$$\sum_{x_{d_1} x_i \in L_2 \setminus \{0\}} (1 + \underline{\lambda}^2 \underline{\mu}^{2i-d} \|x_{d_1} x_i\|_\infty)^{-N_1}.$$

By the integral test, we get the desired estimate for $\Theta_{\mathrm{st},ij}(\Psi, \widetilde{b})$. Note that if $N_1 > 1$,

$$\inf(1, \lambda^{-N_1}) \leq \inf(1, \lambda^{-M})$$

for all $1 \leq M \leq N_1$. Also note that the right hand side of (2.5) does not depend on $x_{d_1}$. □

The proof of the following proposition is similar to that of (2.2) assuming (2.4), and is left to the reader.

**Proposition (2.6)** *For any $N \geq 1$,*

$$\int_\mathbb{A} |\Theta_{\mathrm{st},ij}(\Psi, \widetilde{b})| du \ll \begin{cases} \lambda^{-(j-i+2N)} \mu^{(2i-d)N-(2j-d)+2} & \text{if } \lambda \leq 1, \mu \leq 1 \\ \mathrm{rd}_{2,N}(\lambda, \mu) & \text{otherwise} \end{cases}$$

**Definition (2.7)** *For $\Psi \in \mathscr{S}(\mathbb{A}^{j-d_1+1})$, we define*

$$\Theta_{\mathrm{st},i,1}(\Psi, \widetilde{b}) = \sum_{\substack{x \in Z'^{\mathrm{ss}}_{0k} \\ x_i \in k^\times}} \Psi(A'_{\mathrm{st},i}(\widetilde{b}, x), \underline{\lambda}\underline{\mu}^{d-2i}(x_i + c_i x_{d_1} u^{i-d_1})).$$

Note that $\Theta_{\mathrm{st},i,1}(\Psi, \widetilde{b})$ is a function on $\widetilde{B}_\mathbb{A}/T_k$.
By a similar method as in (2.4), we get the following estimate.

**Proposition (2.8)** *There exists a Schwartz–Bruhat function $0 \leq \Psi' \in \mathscr{S}(\mathbb{A}^2)$ such that for any $N \geq 1$,*

$$\Theta_{\mathrm{st},i,1}(\Psi, \widetilde{b}) \ll \sum_{x_{d_1} \in k^\times} \Psi'(\underline{\lambda} x_{d_1}, \underline{\lambda}\underline{\mu}^{-2} x_{d_1} u)(\lambda^2 \mu^{d-2i})^{-N}.$$

We get the following proposition by Lemma (1.2.6) in [yukiec].

**Proposition (2.9)** *For any $N \geq 1$,*

$$\int_\mathbb{A} |\Theta_{\mathrm{st},i,1}(\Psi, \widetilde{b})| du \ll \lambda^{-(2N+2)} \mu^{(2i-d)N+2}$$

**Definition (2.10)** *For $\Psi \in \mathscr{S}(\mathbb{A}^{j-d_1+1})$, we define*

$$\Theta_{\mathrm{st},i,2}(\Psi, \widetilde{b}) = \sum_{\substack{x \in Z'^{\mathrm{ss}}_{0k} \\ x_i \in k^\times}} \Psi(A'_{\mathrm{st},i}(\widetilde{b}, x), \underline{\lambda}^{-1}\underline{\mu}^{2i-d} x_i) \langle -c_i x_{d_1} x_i u^{i-d_1} \rangle.$$

Note that $\Theta_{\mathrm{st},i,2}(\Psi, \widetilde{b})$ is a function on $\widetilde{B}_\mathbb{A}/T_k$.

**Proposition (2.11)** *There exists a Schwartz–Bruhat function $0 \leq \Psi' \in \mathscr{S}(\mathbb{A}^2)$ such that for any $N \geq 1$,*

$$\Theta_{\mathrm{st},i,2}(\Psi, \widetilde{b}) \ll \sum_{x_{d_1} \in k^\times} \Psi'(\underline{\lambda} x_{d_1}, \underline{\lambda}\underline{\mu}^{-2} x_{d_1} u)(\lambda \mu^{d-2i})^N.$$

*Proof.* There exist Schwartz–Bruhat functions $0 \leq \Psi' \in \mathscr{S}(\mathbb{A}^2)$ and $0 \leq \Psi_1 \in \mathscr{S}(\mathbb{A})$ such that

$$\Psi(A'_{\mathrm{st},i}(\widetilde{b}, x)) \ll \Psi'(\underline{\lambda} x_{d_1}, \underline{\lambda}\underline{\mu}^{-2} x_{d_1} u) \Psi_1(\underline{\lambda}^{-1} \underline{\mu}^{2i-d} x_i).$$

So Proposition (2.11) follows from Lemma (1.2.6) in [yukiec]. $\square$

By (2.11), we get the following proposition.

**Proposition (2.12)** *For any $N_1, N_2 \geq 1$,*

$$\int_\mathbb{A} |\Theta_{\mathrm{st},i,2}(\Psi, \widetilde{b})| du \ll \lambda^{-(N_1-N_2+1)} \mu^{(d-2i)N_2+2}.$$

### §3 Unstable distributions

In this section, we define distributions which will be needed later, and prove their convergence.

First, we recall the definition of the Tate zeta function. For $\Psi \in \mathscr{S}(\mathbb{A})$, $\omega$ is a character of $\mathbb{A}^\times/k^\times$, $t \in \mathbb{A}^\times$, and $s \in \mathbb{C}$, we define

$$\Theta_1(\Psi, t) = \sum_{x \in k^\times} \Psi(tx),$$

$$\Sigma_1(\Psi, \omega, s) = \int_{\mathbb{A}^\times/k^\times} \omega(t)|t|^s d^\times t.$$

Properties of $\Sigma_1(\Psi, \omega, s)$ are well known.

For the rest of this section, $\omega = (\omega_1, \omega_2)$ is a character of $\widetilde{G}_\mathbb{A}/G_k$.

**Definition (3.1)** Let $\Phi \in \mathscr{S}(V_\mathbb{A})$, $x = (x_i, \cdots, x_{j-1}) \in \mathbb{A}^{j-i}$, $x_j \in \mathbb{A}$, $y = (y_j, \cdots, y_d)$, and $dy = dy_j \cdots dy_d$. For $d_1 \leq i < j \leq d$, we define

(1) $$\widetilde{R}_{ij}\Phi(x, x_j) = \int_{\mathbb{A}^{d-j+1}} \Phi(\overbrace{0, \cdots, 0}^{i}, x, y)\langle x_j y_j \rangle dy,$$

(2) $$R_{ij}\Phi(x) = \widetilde{R}_{ij}\Phi(x, 0).$$

For convenience, we define

$$R_{id+1}\Phi(x) = \Phi(\overbrace{0, \cdots, 0}^{i}, x).$$

**Definition (3.2)** Let $[\frac{d}{2}] \leq i < j \leq d$. For $\Psi \in \mathscr{S}(\mathbb{A}^{j-i+1})$, $\lambda \in \mathbb{R}_+$, and a complex variable $s$, we define

$$\Sigma_{ij}(\Psi, \omega, \lambda, s) = \int_{B_\mathbb{A}^0/B_k} \mu^s \omega(b^0) \Theta_{ij}(\Psi, (\underline{\lambda}, b^0)) db^0.$$

By (2.2), we get the following proposition.

**Proposition (3.3)** *For a fixed $\underline{\lambda}$, $\Sigma_{ij}(\Psi, \omega, \lambda, s)$ is an entire function of $s$ and is bounded on any vertical strip.*

Note that $A_{ii+1}(\widetilde{b}, x)$ does not depend on $u$ and $\langle -x_j p_{ii+1}(x, u) \rangle$ is a non-trivial character of $u \in \mathbb{A}/k$. Therefore, $\Sigma_{ii+1}(\Psi, \omega, \lambda, s) = 0$.

For the rest of this section, we assume that $d$ is even, i.e. $d = 2d_1$.

**Definition (3.4)** Let $d_1 + 2 \leq i < j \leq d$. For $\Psi \in \mathscr{S}(\mathbb{A}^{j-d_1+1})$, $\lambda \in \mathbb{R}_+$, and a complex variable $s$, we define

$$\Sigma_{\mathrm{st},ij}(\Psi, \omega, \lambda, s) = \int_{B_\mathbb{A}^0/B_k} \mu^s \omega(b^0) \Theta_{\mathrm{st},ij}(\Psi, (\underline{\lambda}, b^0)) db^0.$$

By (2.6), we get the following proposition.

**Proposition (3.5)** *For a fixed $\underline{\lambda}$, $\Sigma_{\mathrm{st},ij}(\Psi,\omega,\underline{\lambda},s)$ is an entire function of $s$ and is bounded on any vertical strip.*

Let $d = 2d_1$, $d_1 + 2 \leq i \leq d$.

**Definition (3.6)** *Let $d_1 + 2 \leq i \leq d$. For $\Psi \in \mathscr{S}(\mathbb{A}^{i-d_1+1})$, $\lambda \in \mathbb{R}_+$, and a complex variable $s$, we define*

$$\Sigma_{\mathrm{st},i,1}(\Psi,\omega,\lambda,s) = \int_{\substack{B_{\mathbb{A}}^0/B_k \\ \mu \leq 1}} \mu^s \omega(b^0) \Theta_{\mathrm{st},i,1}(\Psi,(\underline{\lambda},b^0)) db^0.$$

By (2.9), we get the following proposition.

**Proposition (3.7)** *For a fixed $\lambda$, $\Sigma_{\mathrm{st},i,1}(\Psi,\omega,\lambda,s)$ is an entire function of $s$ and is bounded on any vertical strip.*

**Definition (3.8)** *Let $d_1 + 2 \leq i \leq d$. For $\Psi \in \mathscr{S}(\mathbb{A}^{i-d_1+1})$, $\lambda \in \mathbb{R}_+$, and a complex variable $s$, we define*

$$\Sigma_{\mathrm{st},i,2}(\Psi,\omega,\lambda,s) = \int_{\substack{B_{\mathbb{A}}^0/B_k \\ \mu \geq 1}} \mu^s \omega(b^0) \Theta_{\mathrm{st},i,2}(\Psi,(\underline{\lambda},b^0)) db^0.$$

By (2.12), we get the following proposition.

**Proposition (3.9)** *For a fixed $\lambda$, $\Sigma_{\mathrm{st},i,2}(\Psi,\omega,\lambda,s)$ is an entire function of $s$ and is bounded on any vertical strip.*

We recall the fundamental domain for $G_{\mathbb{A}}^0/H_k$.

**Definition (3.10)** *For $u \in \mathbb{A}$, we define $\alpha(u) = |t_1({}^t n(u))|$.*

Note that if $u = (u_v)_v \in \mathbb{A}$, $\alpha(u) = \prod_v \alpha_v(u_v)$ where

$$\alpha_v(u_v) = \begin{cases} (1+|u_v|_v^2)^{-\frac{1}{2}} & v \text{ is a real place,} \\ (1+|u_v|_v)^{-1} & v \text{ is an imaginary place,} \\ \sup(1,|u_v|_v^{-1})^{-1} & v \text{ is a finite place.} \end{cases}$$

**Proposition (3.11)** *We define*

$$X_H = \{g^0 \in G_{\mathbb{A}}^0 \mid |t_1(g^0)| \geq \sqrt{\alpha(t_1(g^0)^{-1} t_2(g^0) u(g^0))}\}.$$

We proved in Proposition (2.6) in [yukiea] that for any measurable function $f(g)$ on $G_{\mathbb{A}}^0/H_k$,

$$\int_{G_{\mathbb{A}}^0/H_k} f(g) dg = \int_{X_H/T_k} f(g) dg$$

if the right hand side converges absolutely.

**Definition (3.12)** *For $\Phi \in \mathscr{S}(V_{\mathbb{A}})$ and complex variables $s, s_1$, we define*

(1) $T_V(\Phi, \omega, s, s_1) = \displaystyle\int_{\mathbb{R}_+ \times \widehat{B}_{\mathbb{A}}^0/T_k} \omega(\widehat{t}^{\,0}) \Theta_{Z_0'}(\Psi, d(\lambda)n(u_0)\widehat{t}^{\,0}) \lambda^s \alpha(u_0)^{s_1} d^\times \lambda d\widehat{b}^{\,0}$,

(2) $T_{V+}(\Phi, \omega, s, s_1) = \displaystyle\int_{\substack{\mathbb{R}_+ \times \widehat{B}_{\mathbb{A}}^0/T_k \\ \lambda \geq 1}} \omega(\widehat{t}^{\,0}) \Theta_{Z_0'}(\Psi, d(\lambda)n(u_0)\widehat{t}^{\,0}) \lambda^s \alpha(u_0)^{s_1} d^\times \lambda d\widehat{b}^{\,0}$,

(3) $T_V^1(\Phi, \omega, s_1) = \displaystyle\int_{\widehat{B}_{\mathbb{A}}^0/T_k} \omega(\widehat{t}^{\,0}) \Theta_{Z_0'}(\Psi, n(u_0)\widehat{t}^{\,0}) \alpha(u_0)^{s_1} d^\times \lambda d\widehat{b}^{\,0}$.

By the same argument as in Proposition (2.13) in [yukiea], we get the following proposition.

**Proposition (3.13)** *Let $\epsilon > 0$ be a constant. Then the integral $(3.12)(1)$ converges absolutely and locally uniformly for $\operatorname{Re}(s) \geq 2 + \epsilon$, $\operatorname{Re}(s_1) \geq -\epsilon$, and the integrals $(3.12)(2)$ and $(3.12)(3)$ converge absolutely and locally uniformly for all $s, s_1$.*

**Definition (3.14)** *We define*

$$T_V(\Phi, \omega, s) = \left.\frac{d}{ds_1}\right|_{s_1=0} T_V(\Phi, \omega, s, s_1),$$

$$T_{V+}(\Phi, \omega, s) = \left.\frac{d}{ds_1}\right|_{s_1=0} T_{V+}(\Phi, \omega, s, s_1),$$

$$T_V^1(\Phi, \omega) = \left.\frac{d}{ds_1}\right|_{s_1=0} T_V^1(\Phi, \omega, s_1).$$

## §4 The smoothed Eisenstein series

In this section, we define distributions related to the smoothed Eisenstein series.

Let $z = (z_1, z_2) \in \mathbb{C}^2$, and $z_1 + z_2 = 0$. The point $\rho = (\frac{1}{2}, -\frac{1}{2})$ is known as half the sum of positive weights. The Eisenstein series of $G_{2\mathbb{A}}^0$ for $B_2$ is defined as

$$E(g_2, z) = \sum_{\gamma \in G_k/B_k} a(t_1(g_2\gamma), t_1(g_2\gamma))^{z+\rho}$$

for $g_2 \in G_{2\mathbb{A}}^0$. For the analytic properties of $E(g_2, z)$, see [wright], [yukiec].

We identify the Weyl group of $G$ with the group of permutation matrices in $G_2$, which can be identified with the permutation group of two numbers $\{1, 2\}$ also. If $\tau$ is such a permutation and $z = (z_1, z_2)$ is as above, we define $\tau z = (z_{\tau(1)}, z_{\tau(2)})$. In our situation, $\tau$ is either 1 or the transposition $(12)$. It is proved in Lemma (2.4.13) in [yukiec] that $E(g_2, z) = E(g_2^\iota, z)$ ($g_2^\iota$ is defined in §1). Note that $\nu$ in this paper corresponds to $\tau_G$ in Lemma (2.4.13) [yukiec], and $-\nu z = z$ if $z_1 + z_2 = 0$.

Let $\psi(z)$ be an entire function of $z$ such that

$$\sup_{c_1 < \operatorname{Re}(z) < c_2} (1+|z|)^N |\psi(z)| < \infty$$

for all $c_1 < c_2, N > 0$. For a complex variable $w$, we define

$$\Lambda_\psi(w; z) = \frac{\psi(z)}{w - (z_1 - z_2)}.$$

When there is no confusion we drop $\psi$ and use the notation $\Lambda(w; z)$ instead. Note that $\Lambda(w; \rho) = \frac{\psi(\rho)}{w-1}$.

**Definition (4.1)** For a complex variable $w$, we define

$$\mathscr{E}(g^0, w, \psi) = \frac{1}{2\pi\sqrt{-1}} \int_{1 < \text{Re}(z_1 - z_2) = q < \text{Re}(w)} E(g_2, z)\Lambda_\psi(w; z) d(z_1 - z_2).$$

The function $\mathscr{E}(g^0, w, \psi)$ is called a smoothed Eisenstein series. When there is no confusion we drop $\psi$ and use the notation $\mathscr{E}(g^0, w)$ instead.

**Definition (4.2)** *Let $f(w), g(w)$ be holomorphic functions of $w \in \mathbb{C}$ in some right half plane. We use the notation $f(w) \sim g(w)$ if $f(w) - g(w)$ can be continued meromorphically to $\{w \mid \text{Re}(w) > 1 - \epsilon\}$ for some $\epsilon > 0$ and is holomorphic at $w = 1$.*

Let $\Omega \subset \widehat{B}_{\mathbb{A}}^0$ be a compact set and $\eta > 0$ a constant. A set of the form

$$\mathfrak{S}^0 = \{ka(\underline{\mu}, \underline{\mu}^{-1})n \mid k \in K, \, \mu \geq \eta, \, n \in \Omega\}$$

is called a Siegel set. We choose $\Omega$ large enough and $\eta$ small enough so that $\mathfrak{S}^0$ surjects to $G_{\mathbb{A}}^0/G_k$.

The following proposition was first proved for GL(2) (which is all we need in this paper) by Shintani in [shintania] and later generalized to GL(n) by the author in [yukiec].

**Proposition (4.3)** *(1) If $f(g^0)$ is a function of $g^0 \in G_{\mathbb{A}}^0/G_k$ such that there is a constant $r < 2$ and*

$$f(ka(\underline{\mu}, \underline{\mu}^{-1})n(u)) \ll \mu^r$$

*for $ka(\underline{\mu}, \underline{\mu}^{-1})n(u) \in \mathfrak{S}^0$, the integral*

$$\int_{G_{\mathbb{A}}^0/G_k} f(g)\mathscr{E}(g^0, w)dg^0$$

*becomes a holomorphic function for $\text{Re}(w) \geq 1 - \epsilon$ for a constant $\epsilon > 0$ except possibly for a simple pole at $w = 1$ with residue*

$$\varrho\psi(\rho) \int_{G_{\mathbb{A}}^0/G_k} f(g^0)dg^0.$$

*(2) If $f(g^0)$ is a slowly increasing function of $g^0 \in G_{\mathbb{A}}^0/G_k$, the integral*

$$\int_{G_{\mathbb{A}}^0/G_k} f(g^0)\mathscr{E}(g^0, w)dg^0$$

becomes a holomorphic function on a certain right half plane.

(3) If $\omega_2$ is a character of $\mathbb{A}^\times/k^\times$,
$$\int_{G_\mathbb{A}^0/G_k} \omega(g^0)\mathscr{E}(g^0,w)dg^0 = \delta(\omega_1)\delta(\omega_2)\Lambda(w;\rho)$$
($\omega = (\omega_1,\omega_2)$).

We included the first factor $\mathbb{A}^1$ in the statement instead of considering just $\mathrm{GL}(2)$, but it does not make any difference because $\mathbb{A}^1/k^\times$ is compact and the volume of $\mathbb{A}^1/k^\times$ is 1. Roughly speaking, this proposition says that we can multiply $\mathscr{E}(g,w)$ to any slowly increasing function and make it integrable if $\mathrm{Re}(w) \gg 0$. Moreover, if the function is integrable to begin with, the resulting integral, as a function of $w$, has a simple pole at $w = 1$ with residue a constant multiple of the original integral. This is the key idea to separate contributions from strata.

Let $I(\Phi,\omega,w)$ be the following integral
$$\int_{G_\mathbb{A}^0/G_k} \omega(g^0)\left(\sum_{V_k\setminus V_k^\mathrm{s}} \widehat{\Phi}((g^0)^\iota x) - \sum_{V_k\setminus V_k^\mathrm{s}} \Phi(g^0 x)\right)\mathscr{E}(g^0,w)dg^0.$$

By (4.3),
$$I(\Phi,\omega,w) \sim \varrho\Lambda(w;\rho)I(\Phi,\omega).$$

**Definition (4.4)** *For a complex variable $w$, we define*

(1) $\qquad I_i(\Phi,\omega,w) = \displaystyle\int_{G_\mathbb{A}^0/G_k} \omega(g^0)\Theta_{S_i}(\Phi,g^0)\mathscr{E}(g^0,w)dg^0,$

(2) $\qquad I_{\mathrm{st},i}(\Phi,\omega,w) = \displaystyle\int_{G_\mathbb{A}^0/G_k} \omega(g^0)\Theta_{S_{\mathrm{st},i}}(\Phi,g^0)\mathscr{E}(g^0,w)dg^0,$

(3) $\qquad I_\mathrm{st}(\Phi,\omega,w) = \displaystyle\int_{G_\mathbb{A}^0/G_k} \omega(g^0)\Theta_{V^\mathrm{ss}_{\mathrm{st}k}}(\Phi,g^0)\mathscr{E}(g^0,w)dg^0,$

(4) $\qquad I_\#(\Phi,\omega,w) = \Phi(0)\displaystyle\int_{G_\mathbb{A}^0/G_k} \omega(g^0)\mathscr{E}(g^0,w)dg^0,$

*where $[\frac{d}{2}] < i \leq d$ in (1), and $d$ is even ($d = 2d_1$) and $d_1 + 2 \leq i \leq d$ in (2).*

Since $\Theta_{S_i}(\Phi,g^0)$ etc. are slowly increasing functions, by (4.3), the above integrals converge absolutely if $\mathrm{Re}(w) \gg 0$. Since $E(g_2^\iota,z) = E(g_2,z)$ for $g_2 \in G_{2\mathbb{A}}^0$, $\mathscr{E}(g^0,w) = \mathscr{E}((g^0)^\iota,w)$ for $g^0 \in G_\mathbb{A}^0$. Therefore, by (1.4),

$$\begin{aligned}
(4.5)\quad I(\Phi,\omega,w) &= I_\#(\widehat{\Phi},\omega^{-1},w) - I_\#(\Phi,\omega,w) \\
&\quad + \sum_{[\frac{d}{2}]<i\leq d}(I_i(\widehat{\Phi},\omega^{-1},w) - I_i(\Phi,\omega,w)) \text{ if } d \text{ is odd,} \\
I(\Phi,\omega,w) &= I_\#(\widehat{\Phi},\omega^{-1},w) - I_\#(\Phi,\omega,w) \\
&\quad + \sum_{d_1<i\leq d}(I_i(\widehat{\Phi},\omega^{-1},w) - I_i(\Phi,\omega,w)) \\
&\quad + I_\mathrm{st}(\widehat{\Phi},\omega^{-1},w) - I_\mathrm{st}(\Phi,\omega,w)) \\
&\quad + \sum_{d_1+2\leq i\leq d}(I_{\mathrm{st},i}(\widehat{\Phi},\omega^{-1},w) - I_{\mathrm{st},i}(\Phi,\omega,w)) \text{ if } d = 2d_1.
\end{aligned}$$

By (4.3)(3),

(4.6)    $$I_\#(\Phi,\omega,w) = \delta(\omega_1)\delta(\omega_2)\Phi(0)\Lambda(w;\rho).$$

We study each term of the right hand side of (4.5) in §§5,6. For the rest of this section, we consider estimates of the smoothed Eisenstein series, and introduce notations similar to those in §3.6 in [yukiec].

The function $\mathscr{E}(g^0, w)$ is defined on $G_\mathbb{A}^0/G_k$. So we can consider its Fourier expansion. We denote the constant term and the non-constant term of the Fourier expansion of $\mathscr{E}(g^0, w)$ by $\mathscr{E}_0(g^0, w)$ and $\mathscr{E}_1(g^0, w)$ respectively. The following lemma is a special case of Proposition (3.4.30)(1) in [yukiec] (it also follows easily from Proposition (2.3.26) [yukiec]).

**Proposition (4.7)** *Let $\delta > 0$. Then for any $l > 1$,*

$$\mathscr{E}_1(g^0, w) \ll |t_1(g^0)t_2(g^0)^{-1}|^{\frac{1}{2}-l}$$

*for $g^0 \in G_\mathbb{A}^0$, $\mathrm{Re}(w) \geq 1 - \delta$.*

In Proposition (3.4.30)(1) [yukiec], we only have the part $t(g^0)^{\nu(\tau q + \rho - s_I(l))}$ and $\nu = 1, \tau = (12)$, and $s_I(l)$ consists of a single number $l > 1$. Then we choose $\tau q$ close enough to the origin. Note that $\rho$ contributes $\frac{1}{2}$ in the above proposition. In this paper, we don't exactly need the above precise bound $l > 1$. Since the formulation in [yukiec] is quite general, the reader may see [wright] for a similar estimate also. Even though the statement in [wright] is not the optimum one, a careful application of the argument in [wright] implies the bound in the above proposition.

The Weyl group of $G$ is $\mathfrak{S}_2$ (the permutation group of $\{1, 2\}$). Let $\tau$ be a Weyl group element, i.e. $\tau = 1$ or $(12)$. Let $s_\tau = z_2 - z_1$ if $\tau = 1$, and $s_\tau = z_1 - z_2$ if $\tau = (12)$. Let $\phi(s)$ be the function we defined in §1. Let $M_\tau(s_\tau) = 1$ if $\tau = 1$, and $M_\tau(s_\tau) = \phi(s_\tau)$ if $\tau = (12)$. For a fixed $\tau$, any function of $z$ can be considered as a function of $s_\tau$. So we define $\widetilde{\Lambda}_\tau(w; s_\tau) = M_\tau(s_\tau)\Lambda(w;z)$, where we consider $\Lambda(w;z)$ as a function of $s_\tau$. Explicitly,

$$\widetilde{\Lambda}_\tau(w; s_\tau) = \begin{cases} \frac{\psi(-\frac{1}{2}s_\tau, \frac{1}{2}s_\tau)}{w+s_\tau} & \tau = 1, \\ \phi(s_\tau)\frac{\psi(\frac{1}{2}s_\tau, -\frac{1}{2}s_\tau)}{w-s_\tau} & \tau = (12). \end{cases}$$

Suppose $g^0 = (t, g_2)$, where $g_2 = ka(\underline{\mu}t_1, \underline{\mu}^{-1}t_2)n(u)$ is the Iwasawa decomposition and $t, t_1, t_2 \in \mathbb{A}^1$. Then

$$\mathscr{E}_0(g^0, w) = \sum_{\tau=1,(12)} \frac{1}{2\pi\sqrt{-1}} \int_{\mathrm{Re}(w) > \mathrm{Re}(s_\tau) = r > 1} \mu^{1-s_\tau} \widetilde{\Lambda}_\tau(w; s_\tau) ds_\tau.$$

This is the reason why we introduced the above somewhat general formulation. Even if we are dealing with the simplest possible non-trivial group $\mathrm{GL}(2)$, this formulation has an advantage that the right hand side of the above formula can be treated simultaneously for $\tau = 1, (12)$. Since we are going to deal with many terms, without this formulation, we may have to write two terms each for contributions from the constant term $\mathscr{E}_0(g^0, w)$.

For a complex variable $s$, we define

$$\widetilde{\Lambda}(w; s) = \frac{\psi(\frac{1}{2}s, -\frac{1}{2}s)}{w - s}.$$

**Proposition (4.8)** *For any $\epsilon > 0$, there exists $\delta > 0$ and $|c_1|, |c_2| < \epsilon$ such that for any $l > 1$,*

$$|\mathscr{E}(g^0, w) - \varrho\Lambda(w; \rho)| \ll |t_1(g^0)|^{c_1}|t_2(g^0)|^{c_2} + |t_1(g^0)t_2(g^0)^{-1}|^{\frac{1}{2}-l}$$

*for $g \in G_{\mathbb{A}}^0$, $\mathrm{Re}(w) \geq 1 - \delta$.*

*Proof.* By (4.7), we only have to consider the difference between $\mathscr{E}_0(g^0, w)$ and $\varrho\Lambda(w; \rho)$. If $\tau = 1$, we can choose $r$ close to 1 and $\widetilde{\Lambda}_\tau(w; s_\tau)$ is holomorphic at $w = 1$. If $\tau = (12)$,

$$\frac{1}{2\pi\sqrt{-1}} \int_{\mathrm{Re}(w) > \mathrm{Re}(s_\tau) = r > 1} \mu^{1-s_\tau} \widetilde{\Lambda}_\tau(w; s_\tau) ds_\tau$$
$$= \varrho\Lambda(w; \rho) + \frac{1}{2\pi\sqrt{-1}} \int_{\mathrm{Re}(w) > \mathrm{Re}(s_\tau) = r < 1} \mu^{1-s_\tau} \widetilde{\Lambda}_\tau(w; s_\tau) ds_\tau.$$

Then we can choose $r < 1$ close to 1. $\square$

In §§5,6, we consider contour integrals of the form $\int_{\mathrm{Re}(s_\tau) = r} \cdots \widetilde{\Lambda}_\tau(w; s_\tau) ds_\tau$. We assume that whenever we consider such an integral, we only consider $w$ such that $\mathrm{Re}(w - (z_1 - z_2)) > 0$. This ensures that the denominator in $\Lambda(w; z)$ is non-zero. This condition is satisfied if $\mathrm{Re}(w) > -r$ (resp. $\mathrm{Re}(w) > r$) if $\tau = 1$ (resp. $\tau = (12)$).

## §5 Contributions from unstable strata

In this section, we consider each term in (4.5).

We define

(5.1) $\quad \delta_i(\omega) = \delta(\omega_2 \omega_1^{i-d})\delta(\omega_2 \omega_1^{-i}),$

$$\Sigma_{S_i}(\Phi, \omega, s_\tau) = \frac{\delta_i(\omega)}{2i - d} \Sigma_1\left(R_{ii+1}\Phi, \omega_1, \frac{-s_\tau - 1 + (i+1)(d-i)}{d - 2i}\right),$$

$$\Sigma_{S_i}(\Phi, \omega) = \frac{\delta_i(\omega)}{2i - d} \Sigma_1\left(R_{ii+1}\Phi, \omega_1, \frac{(i+1)(d-i) - 2}{d - 2i}\right).$$

**Proposition (5.2)** *For $[\frac{d}{2}] < i \leq d$,*

$$I_i(\Phi, \omega, w) \sim \sum_{i+1 < j \leq d} \Sigma_{ij}(\widetilde{R}_{ij}\Phi, \omega, 1, j(d - j + 1))\varrho\Lambda(w; \rho)$$
$$+ \Sigma_{S_i}(\Phi, \omega)\varrho\Lambda(w; \rho).$$

*Proof.*
$$I_i(\Phi, \omega, w) = \int_{G^0_{\mathbb{A}}/G_k} \omega(g^0) \Theta_{S_i}(\Phi, g^0) \mathscr{E}(g^0, w) dg^0$$
$$= \int_{G^0_{\mathbb{A}}/B_k} \omega(g^0) \Theta_{Y_i}(\Phi, g^0) \mathscr{E}(g^0, w) dg^0$$
$$= \int_{B^0_{\mathbb{A}}/B_k} \omega(b^0) \Theta_{Y_i}(\Phi, b^0) \mathscr{E}(b^0, w) db^0.$$

Note that
$$\Theta_{Y_i}(\Phi, b^0) = \sum_{\substack{x \in Y'_{idk} \\ x_d \in k}} \Phi(A_{id}(b^0, x), \gamma_d(t^0)(x_d + p_{id}(x, u))).$$

We will apply the Poisson summation formula with respect to the variable $x_d$. Fixing $x_i, \cdots, x_{d-1}$ and $b^0$, the partial Fourier transform of $\Phi(A_{id}(b^0, x), \gamma_d(t^0)(x_d + p_{id}(x, u)))$ with respect to $x_d$ is

$$\int_{\mathbb{A}} \Phi(A_{id}(b^0, x), \gamma_d(t^0)(y + p_{id}(x, u))) \langle x_d y \rangle dy$$
$$= \langle -x_d p_{id}(x, u) \rangle \int_{\mathbb{A}} \Phi(A_{id}(b^0, x), \gamma_d(t^0) y) \langle x_d y \rangle dy$$
$$= \langle -x_d p_{id}(x, u) \rangle \mu^{-d} \widetilde{R}_{id} \Phi(A_{id}(b^0, x), \gamma_d(t^0)^{-1} x_d).$$

Therefore, applying the Poisson summation formula to the variable $x_d$,
$$\Theta_{Y_i}(\Phi, b^0) = \sum_{\substack{x \in Y'_{idk} \\ x_d \in k}} \langle -x_d p_{id}(x, u) \rangle \mu^{-d} \widetilde{R}_{id} \Phi(A_{id}(b^0, x), \gamma_d(t^0)^{-1} x_d)$$
$$= \mu^{-d} \Theta_{id}(\widetilde{R}_{ij}\Phi, b^0) + \sum_{x \in Y'_{idk}} \mu^{-d} R_{id} \Phi(A_{id}(b^0, x)).$$

Continuing this process,
$$\Theta_{Y_i}(\Phi, b^0) = \sum_{i < j \leq d} \mu^{j(d-j+1)} \Theta_{ij}(\widetilde{R}_{ij}\Phi, b^0) + \mu^{(i+1)(d-i)} \Theta_1(R_{ii+1}\Phi, \gamma_i(t^0)).$$

(If $i = d$, we do not have the first term.) Note that
$$\mu^{j(d-j+1)} \mu^{2(j-1)-d} = \mu^{(j-1)(d-(j-1)+1)}.$$

By (2.2) and (4.8),
$$\int_{B^0_{\mathbb{A}}/B_k} \omega(b^0) \mu^{j(d-j+1)} \Theta_{ij}(\widetilde{R}_{ij}\Phi, b^0) \mathscr{E}(b^0, w) db^0$$
$$\sim \int_{B^0_{\mathbb{A}}/B_k} \omega(b^0) \mu^{j(d-j+1)} \Theta_{ij}(\widetilde{R}_{ij}\Phi, b^0) db^0 \varrho\Lambda(w; \rho)$$
$$= \Sigma_{ij}(\widetilde{R}_{ij}\Phi, \omega, 1, j(d-j+1)) \varrho\Lambda(w; \rho).$$

Note that $\Sigma_{ii+1}(\widetilde{R}_{ii+1}\Phi,\omega,1,j(d-j+1)) = 0$.

Let $\tau$ be a Weyl group element. Since $\gamma_i(t^0) = \underline{\mu}^{d-2i}tt_1^{d-i}t_2^i$ and $\omega(t^0) = \omega_1(t)\omega_2(t_1t_2)$,

$$\int_{B_\mathbb{A}^0/B_k} \omega(b^0)\mu^{-s_\tau+1+(i+1)(d-i)}\Theta_1(R_{ii+1}\Phi,\gamma_i(t^0))db^0$$
$$= \int_{T_\mathbb{A}^0/T_k} \omega(t^0)\mu^{-s_\tau-1+(i+1)(d-i)}\Theta_1(R_{ii+1}\Phi,\gamma_i(t^0))dt^0$$
$$= \Sigma_{S_i}(\Phi,\omega,s_\tau)$$

if $\operatorname{Re}(s_\tau) \gg 0$.

Therefore,

$$\int_{B_\mathbb{A}^0/B_k} \omega(b^0)\mu^{(i+1)(d-i)}\Theta_1(R_{ii+1}\Phi,\gamma_i(t^0))\mathscr{E}(b^0,w)db^0$$
$$= \sum_{\tau=1,(12)} \frac{1}{2\pi\sqrt{-1}}\int_{\operatorname{Re}(s_\tau)=r\gg 0} \Sigma_{S_i}(\Phi,\omega,s_\tau)\widetilde{\Lambda}_\tau(w;s_\tau)ds_\tau.$$

If $\tau = 1$,

$$\frac{1}{2\pi\sqrt{-1}}\int_{\operatorname{Re}(s_\tau)=r} \Sigma_{S_i}(\Phi,\omega,s_\tau)\widetilde{\Lambda}_\tau(w;s_\tau)ds_\tau \sim 0.$$

Suppose $\tau = (12)$. It is easy to see that $\frac{(i+1)(d-i)-2}{d-2i} \neq 0,1$ for $[\frac{d}{2}] < i \leq d$ (we are still assuming $d \geq 4$). Therefore, by moving the contour to the left as usual,

$$\frac{1}{2\pi\sqrt{-1}}\int_{\operatorname{Re}(s_\tau)=r} \Sigma_{S_i}(\Phi,\omega,s_\tau)\widetilde{\Lambda}_\tau(w;s_\tau)ds_\tau \sim \Sigma_{S_i}(\Phi,\omega)\varrho\Lambda(w;\rho).$$

□

For the rest of this section, we assume $d = 2d_1$.

**Proposition (5.3)** *For $d_1 + 2 \leq i \leq d$,*

$$I_{\mathrm{st},i}(\Phi,\omega,w) \sim \sum_{i<j\leq d} \Sigma_{\mathrm{st},ij}(\widetilde{R}_{ij}\Phi,\omega,1,j(d-j+1))\varrho\Lambda(w;\rho)$$
$$+ \int_{B_\mathbb{A}^0/T_k} \omega(b^0)\mu^{(i+1)(d-i)}\Theta_{\mathrm{st},i,1}(R_{d_1i+1}\Phi,b^0)\mathscr{E}(b^0,w)db^0.$$

*Proof.*

$$I_{\mathrm{st},i}(\Phi,\omega,w) = \int_{G_\mathbb{A}^0/G_k} \omega(g^0)\Theta_{S_{\mathrm{st},i}}(\Phi,g^0)\mathscr{E}(g^0,w)dg^0$$
$$= \int_{G_\mathbb{A}^0/T_k} \omega(g^0)\Theta_{Y_{\mathrm{st},i}}(\Phi,g^0)\mathscr{E}(g^0,w)dg^0$$
$$= \int_{B_\mathbb{A}^0/T_k} \omega(b^0)\Theta_{Y_{\mathrm{st},i}}(\Phi,b^0)\mathscr{E}(b^0,w)db^0.$$

Applying the Poisson summation formula to $x_d, \cdots, x_{i+1}$ successively,

$$\Theta_{Y_{\text{st},i}}(\Phi, b^0) = \sum_{i<j\leq d} \mu^{j(d-j+1)} \Theta_{\text{st},ij}(\widetilde{R}_{d_1 j}\Phi, b^0) + \mu^{(i+1)(d-i)} \Theta_{\text{st},i,1}(R_{d_1 i+1}\Phi, b^0).$$

By (2.6) and (4.8),

$$\int_{B_{\mathbb{A}}^0/T_k} \omega(b^0) \mu^{j(d-j+1)} \Theta_{\text{st},ij}(\widetilde{R}_{ij}\Phi, b^0) \mathscr{E}(b^0, w) db^0$$
$$\sim \int_{B_{\mathbb{A}}^0/T_k} \omega(b^0) \mu^{j(d-j+1)} \Theta_{\text{st},ij}(\widetilde{R}_{ij}\Phi, b^0) db^0 \varrho\Lambda(w; \rho)$$
$$= \Sigma_{\text{st},ij}(\widetilde{R}_{ij}\Phi, \omega, 1, j(d-j+1)) \varrho\Lambda(w; \rho).$$

This proves Proposition (5.3).  $\square$

**Proposition (5.4)**

$$\sum_{i=d_1+2}^{d} \int_{B_{\mathbb{A}}^0/T_k} \omega(b^0) \mu^{(i+1)(d-i)} \Theta_{\text{st},i,1}(R_{d_1 i+1}\Phi, b^0) \mathscr{E}(b^0, w) db^0$$
$$\sim \sum_{i=d_1+2}^{d} \Sigma_{\text{st},i,1}(R_{d_1 i+1}\Phi, \omega, 1, (i+1)(d-i)) \varrho\Lambda(w; \rho)$$
$$+ \sum_{i=d_1+2}^{d} \Sigma_{\text{st},i,2}(\widetilde{R}_{d_1 i}\Phi, \omega, 1, i(d-i+1)) \varrho\Lambda(w; \rho)$$
$$- \frac{\delta(\omega_2 \omega_1^{-d_1})}{(d_1+2)(d_1-1)} \int_{\mathbb{A}^1/k^{\times}} \omega_1(t) \Theta_1(R_{d_1 d_1+1}\Phi, t) d^{\times}t \varrho\Lambda(w; \rho)$$
$$- \frac{1}{2\pi\sqrt{-1}} \int_{\text{Re}(s)=1+\delta} \frac{T_V^1(\Phi, \omega, 0)}{s-1} \phi(s) \widetilde{\Lambda}(w; s) ds,$$

where $\delta > 0$ is a small number.

*Proof.* We divide the integral into two parts

$$\{b^0 \in B_{\mathbb{A}}^0/T_k \mid \mu \leq 1\}, \{b^0 \in B_{\mathbb{A}}^0/T_k \mid \mu \geq 1\}.$$

By (2.9) and (4.8),

$$\int_{\substack{B_{\mathbb{A}}^0/T_k \\ \mu \leq 1}} \omega(b^0) \mu^{(i+1)(d-i)} \Theta_{\text{st},i,1}(R_{d_1 i+1}\Phi, b^0) \mathscr{E}(b^0, w) db^0$$
$$\sim \int_{\substack{B_{\mathbb{A}}^0/T_k \\ \mu \leq 1}} \omega(b^0) \mu^{(i+1)(d-i)} \Theta_{\text{st},i,1}(R_{d_1 i+1}\Phi, b^0) db^0 \varrho\Lambda(w; \rho)$$
$$= \Sigma_{\text{st},i,1}(R_{d_1 i+1}\Phi, \omega, 1, (i+1)(d-i)) \varrho\Lambda(w; \rho).$$

Applying the Poisson summation formula to $x_i$,

$$\Theta_{\mathrm{st},i,1}(R_{d_1 i+1}\Phi, b^0) = \mu^{2i-d}\Theta_{\mathrm{st},i,2}(\widetilde{R}_{d_1 i}\Phi, b^0)$$
$$+ \mu^{2i-d}\sum_{x\in Z'^{\mathrm{ss}}_{0k}} R_{d_1 i}\Phi(A'_{\mathrm{st},i}(b^0, x))$$
$$- \sum_{x\in Z'^{\mathrm{ss}}_{0k}} R_{d_1 i+1}\Phi(A'_{\mathrm{st},i+1}(b^0, x)).$$

Since $\mu^{(i+1)(d-i)}\mu^{2i-d} = \mu^{i(d-i+1)}$,

$$\sum_{i=d_1+2}^{d} \mu^{(i+1)(d-i)}\Theta_{\mathrm{st},i,1}(R_{d_1 i+1}\Phi, b^0)$$
$$= \sum_{i=d_1+2}^{d} \mu^{i(d-i+1)}\Theta_{\mathrm{st},i,2}(\widetilde{R}_{d_1 i}\Phi, b^0)$$
$$+ \mu^{(d_1+2)(d_1-1)}\sum_{x\in Z'^{\mathrm{ss}}_{0k}} R_{d_1 d_1+2}\Phi(A'_{\mathrm{st},d_1+2}(b^0, x)) - \Theta_{Z'_0}(\Phi, b^0).$$

By (2.12) and (4.8),

$$\int_{\substack{B^0_{\mathbb{A}}/T_k \\ \mu\geq 1}} \omega(b^0)\mu^{i(d-i+1)}\Theta_{\mathrm{st},i,2}(\Phi, b^0)\mathscr{E}(b^0, w)db^0$$
$$\sim \Sigma_{\mathrm{st},i,2}(\widetilde{R}_{d_1 i+1}\Phi, \omega, 1, i(d-i+1))\varrho\Lambda(w;\rho).$$

Therefore,

$$\sum_{i=d_1+2}^{d} \int_{\substack{B^0_{\mathbb{A}}/T_k \\ \mu\geq 1}} \omega(b^0)\mu^{(i+1)(d-i)}\Theta_{\mathrm{st},i,1}(\Phi, b^0)\mathscr{E}(b^0, w)db^0$$
$$\sim \sum_{i=d_1+2}^{d} \Sigma_{\mathrm{st},i,2}(\widetilde{R}_{d_1 i+1}\Phi, \omega, 1, i(d-i+1))\varrho\Lambda(w;\rho)$$
$$+ \int_{\substack{B^0_{\mathbb{A}}/T_k \\ \mu\geq 1}} \omega(b^0)\mu^{(d_1+2)(d_1-1)}\sum_{x\in Z'^{\mathrm{ss}}_{0k}} R_{d_1 d_1+2}\Phi(A'_{\mathrm{st},d_1+2}(b^0, x))\mathscr{E}(b^0, w)db^0$$
$$- \int_{\substack{B^0_{\mathbb{A}}/T_k \\ \mu\geq 1}} \omega(b^0)\Theta_{Z'_0}(\Phi, b^0)\mathscr{E}(b^0, w)db^0.$$

**Lemma (5.5)**

(1) $\int_{\substack{B^0_{\mathbb{A}}/T_k \\ \mu\geq 1}} \omega(b^0)\mu^{(d_1+2)(d_1-1)}\sum_{x\in Z'^{\mathrm{ss}}_{0k}} R_{d_1 d_1+2}\Phi(A'_{\mathrm{st},d_1+2}(b^0, x))\mathscr{E}(b^0, w)db^0$

$\sim \int_{\substack{B^0_{\mathbb{A}}/T_k \\ \mu\geq 1}} \omega(b^0)\mu^{(d_1+2)(d_1-1)}\sum_{x\in Z'^{\mathrm{ss}}_{0k}} R_{d_1 d_1+2}\Phi(A'_{\mathrm{st},d_1+2}(b^0, x))\mathscr{E}_0(b^0, w)db^0.$

(2) $\int_{\substack{B^0_{\mathbb{A}}/T_k \\ \mu\geq 1}} \omega(b^0)\Theta_{Z'_0}(\Phi, b^0)\mathscr{E}(b^0, w)db^0 \sim \int_{\substack{B^0_{\mathbb{A}}/T_k \\ \mu\geq 1}} \omega(b^0)\Theta_{Z'_0}(\Phi, b^0)\mathscr{E}_0(b^0, w)db^0.$

*Proof.* There exists a Schwartz–Bruhat function $0 \leq \Psi \in \mathscr{S}(\mathbb{A}^2)$ such that

$$\Theta_{Z'_0}(\Phi, b^0) \ll \sum_{x_{d_1} \in k^\times} \Psi(x_{d_1}, \underline{\mu}^{-2} x_{d_1} u).$$

So

$$\int_\mathbb{A} |\Theta_{Z'_0}(\Phi, b^0)| du, \quad \int_\mathbb{A} \sum_{x \in Z'^{ss}_{0k}} |R_{d_1 d_1 + 2}\Phi(A'_{\text{st},d_1+2}(b^0, x))| du$$

are bounded by a constant multiple of $\mu^2$.

Therefore, by (4.7), there exists $\delta > 0$ such that for any $l > 1$,

$$\int_{\substack{B^0_{\mathbb{A}}/T_k \\ \mu \geq 1}} \mu^{(d_1+2)(d_1-1)} \sum_{x \in Z'^{ss}_{0k}} |R_{d_1 d_1 + 2}\Phi(A'_{\text{st},d_1+2}(b^0, x))\mathscr{E}_1(b^0, w)| db^0$$

$$\ll \int_1^\infty \mu^{(d_1+2)(d_1-1)+3-2l} d^\times \mu < \infty,$$

$$\int_{\substack{B^0_{\mathbb{A}}/T_k \\ \mu \geq 1}} |\Theta_{Z'_0}(\Phi, b^0)\mathscr{E}(b^0, w)| db^0 \ll \int_1^\infty \mu^{3-2l} d^\times \mu < \infty$$

respectively for $\text{Re}(w) \geq 1 - \delta$. This proves Lemma (5.5). $\square$

It is easy to see that

$$\int_\mathbb{A} \sum_{x \in Z'^{ss}_{0k}} R_{d_1 d_1 + 2}\Phi(A'_{\text{st},d_1+2}(b^0, x)) du = \mu^2 \Theta_1(R_{d_1 d_1 + 1}\Phi, tt_1^{d_1} t_2^{d_1}).$$

So

$$\int_{\substack{B^0_{\mathbb{A}}/T_k \\ \mu \geq 1}} \omega(b^0) \mu^{(d_1+2)(d_1-1)} \sum_{x \in Z'^{ss}_{0k}} R_{d_1 d_1 + 2}\Phi(A'_{\text{st},d_1+2}(b^0, x))\mathscr{E}_0(b^0, w) db^0$$

$$= \int_{\substack{T^0_{\mathbb{A}}/T_k \\ \mu \geq 1}} \omega(t^0) \mu^{(d_1+2)(d_1-1)} \Theta_1(R_{d_1 d_1 + 1}\Phi, tt_1^{d_1} t_2^{d_1})\mathscr{E}_0(t^0, w) d^\times t^0$$

$$= \delta(\omega_2 \omega_1^{-d_1}) \int_{\mathbb{A}^1/k^\times} \omega_1(t)\Theta_1(R_{d_1 d_1 + 1}\Phi, t) d^\times t$$

$$\times \sum_{\tau=1,(12)} \frac{1}{2\pi\sqrt{-1}} \int_{\text{Re}(s_\tau)=r \gg 0} \frac{\widetilde{\Lambda}_\tau(w; s_\tau)}{s_\tau - 1 - (d_1 + 2)(d_1 - 1)} ds_\tau.$$

If $\tau = 1$,

$$\int_{\text{Re}(s_\tau)=r \gg 0} \frac{\widetilde{\Lambda}_\tau(w; s_\tau)}{s_\tau - 1 - (d_1 + 2)(d_1 - 1)} ds_\tau \sim 0.$$

Suppose $\tau = (12)$. Since $(d_1 + 2)(d_1 - 1) \neq 0$, by moving the contour to the left as usual,

$$\int_{\text{Re}(s_\tau)=r} \frac{\widetilde{\Lambda}_\tau(w; s_\tau)}{s_\tau - 1 - (d_1 + 2)(d_1 - 1)} ds_\tau \sim -\frac{\varrho \Lambda(w; \rho)}{(d_1 + 2)(d_1 - 1)}.$$

We consider $b^0 = n(u_0)t^0 \in B_{\mathbb{A}}^0$. It is easy to see that $\Theta_{Z_0'}(\Phi, b^0)$ does not depend on $\mu$. Since

$$\int_1^\infty \mathscr{E}_0(b^0, w) d^\times \mu = \sum_{\tau=1,(12)} \int_{\mathrm{Re}(s_\tau)=r>1} \frac{\widetilde{\Lambda}_\tau(w; s_\tau)}{s_\tau - 1} ds_\tau,$$

$$\int_{\substack{B_{\mathbb{A}}^0/T_k \\ \mu \geq 1}} \omega(b^0) \Theta_{Z_0'}(\Phi, b^0) \mathscr{E}_0(b^0, w) db^0$$

$$= \sum_{\tau=1,(12)} \frac{1}{2\pi\sqrt{-1}} \int_{\mathrm{Re}(s_\tau)=r>1} \frac{T_V^1(\Phi, \omega, 0)}{s_\tau - 1} \widetilde{\Lambda}_\tau(w; s_\tau) ds_\tau.$$

If $\tau = 1$,

$$\frac{1}{2\pi\sqrt{-1}} \int_{\mathrm{Re}(s_\tau)=r>1} \frac{T_V^1(\Phi, \omega, 0)}{s_\tau - 1} \widetilde{\Lambda}_\tau(w; s_\tau) ds_\tau \sim 0.$$

If $\tau = (12)$, we choose $r = 1 + \delta > 0$ where $\delta > 0$ is any small constant. This proves Proposition (5.4). $\square$

**Proposition (5.6)**

$$\int_{\mathbb{A}^1/k^\times} \omega_1(t) \Theta_1(R_{d_1 d_1+1} \Phi, t) d^\times t - \int_{\mathbb{A}^1/k^\times} \omega_1^{-1}(t) \Theta_1(R_{d_1 d_1+1} \widehat{\Phi}, t) d^\times t$$

$$= \delta(\omega_1)(R_{d_1 d_1+1}\widehat{\Phi}(0) - R_{d_1 d_1+1}\Phi(0)).$$

*Proof.* Since $M_\omega \Phi = \Phi$, $R_{d_1 d_1+1}\widehat{\Phi}$ is the Fourier transform of $R_{d_1 d_1+1}\Phi$ with respect to the character $\langle\ \rangle$. Therefore, by the Poisson summation formula,

$$\omega(t)(\Theta_1(R_{d_1 d_1+1}\Phi, t) - \Theta_1(R_{d_1 d_1+1}\widehat{\Phi}, t^{-1})) = \omega(t)(R_{d_1 d_1+1}\widehat{\Phi}(0) - R_{d_1 d_1+1}\Phi(0)).$$

Integrating over $\mathbb{A}^1/k^\times$, we get Proposition (5.6). $\square$

**Proposition (5.7)**

$$\int_{G_{\mathbb{A}}^0/G_k} \omega(g^0) \Theta_{V,\mathrm{st}}(\Phi, g^0) dg^0 \sim \frac{1}{2\pi\sqrt{-1}} \int_{\mathrm{Re}(s)=1+\delta} \frac{T_V^1(\Phi, \omega, \frac{1-s}{2})}{s - 1} \phi(s) \widetilde{\Lambda}(w; s) ds,$$

where $\delta > 0$ is a small number.

*Proof.* We consider $b^0 = n(u_0)t^0 \in B_{\mathbb{A}}^0$.

$$\int_{G_{\mathbb{A}}^0/G_k} \omega(g^0) \Theta_{V,\mathrm{st}}(\Phi, g^0) \mathscr{E}(g^0, w) dg^0$$

$$= \int_{G_{\mathbb{A}}^0/H_k} \omega(g^0) \Theta_{Z_0'}(\Phi, g^0) \mathscr{E}(g^0, w) dg^0$$

$$= \int_{X_H/T_k} \omega(b^0) \Theta_{Z_0'}(\Phi, b^0) \mathscr{E}(b^0, w) db^0.$$

We show that we can ignore the non-constant term of $\mathscr{E}(b^0, w)$.

There exists a Schwartz–Bruhat function $0 \leq \Psi \in \mathscr{S}(\mathbb{A}^2)$ such that

$$\Theta_{Z_0'}(\Phi, b^0) \ll \sum_{x_{d_1} \in k^\times} \Psi(x_{d_1}, x_{d_1} u_0).$$

So by (4.7), there exists $\delta > 0$ such that for any $l > 1$

$$\int_{X_H/T_k} |\Theta_{Z_0'}(\Phi, b^0) \mathscr{E}_1(b^0, w)| db^0 \ll \int_{\mathbb{A}^1 \times \mathbb{A}} \Psi(q, qu_0) \frac{\alpha(u_0)^{\frac{1-2l}{2}}}{2l - 1} d^\times q du_0$$

for $\operatorname{Re}(w) \geq 1 - \delta$. Since this integral converges absolutely by Proposition (2.13) in [yukiea],

$$\int_{X_H/T_k} \Theta_{Z_0'}(\Phi, b^0) \mathscr{E}(b^0, w) db^0 \sim \int_{X_H/T_k} \Theta_{Z_0'}(\Phi, b^0) \mathscr{E}_0(b^0, w) db^0.$$

Since

$$\int_{\mu \geq \sqrt{\alpha(u_0)}} \mu^{1-s_\tau} d^\times \mu = \frac{\alpha(u_0)^{\frac{1-s_\tau}{2}}}{s_\tau - 1},$$

$$\int_{\widehat{B^0}/T_k} \omega(\widehat{b}^0) \Theta_{Z_0'}(\Phi, \widehat{b}^0) \alpha(u_0)^{\frac{1-s_\tau}{2}} d\widehat{b}^0 = T_V^1\left(\Phi, \omega, \frac{1 - s_\tau}{2}\right),$$

the above integral is equal to

$$\sum_{\tau = 1, (12)} \frac{1}{2\pi\sqrt{-1}} \int_{\operatorname{Re}(s_\tau) = r > 1} \frac{T_V^1(\Phi, \omega, \frac{1-s_\tau}{2})}{s_\tau - 1} \widetilde{\Lambda}_\tau(w; s_\tau) ds_\tau.$$

If $\tau = 1$,

$$\frac{1}{2\pi\sqrt{-1}} \int_{\operatorname{Re}(s_\tau) = r > 1} \frac{T_V^1(\Phi, \omega, \frac{1-s_\tau}{2})}{s_\tau - 1} \widetilde{\Lambda}_\tau(w; s_\tau) ds_\tau \sim 0.$$

If $\tau = (12)$, we choose $r = 1 + \delta$ where $\delta > 0$ is any small constant. This proves Proposition (5.7). $\square$

Since
$$\frac{T^1(\Phi, \omega, \frac{1-s}{2}) - T^1(\Phi, \omega, 0)}{s - 1}$$
is holomorphic at $s = 1$ and the value at $s = 1$ is $-\frac{1}{2} T^1(\Phi, \omega)$, we get the following proposition.

**Proposition (5.8)** *Let $\delta > 0$ be a constant. Then*

$$\frac{1}{2\pi\sqrt{-1}} \int_{\operatorname{Re}(s) = 1 + \delta} \frac{T^1(\Phi, \omega, \frac{1-s}{2}) - T^1(\Phi, \omega, 0)}{s - 1} \phi(s) \widetilde{\Lambda}(w; s) ds$$

$$\sim -\frac{1}{2} T^1(\Phi, \omega) \varrho \Lambda(w; \rho).$$

We define

$$F_1(\Phi,\omega) = \sum_{[\frac{d}{2}]<i\leq d} \Sigma_{S_i}(\Phi,\omega),$$

$$F_2(\Phi,\omega) = \sum_{[\frac{d}{2}]<i\leq d} \sum_{i+2\leq j\leq d} \Sigma_{ij}(\widetilde{R}_{ij}\Phi,\omega,1,j(d-j+1)).$$

If $d = 2d_1$, we define

$$F_3(\Phi,\omega) = \sum_{d_1+2\leq i<j\leq d} \Sigma_{\mathrm{st},ij}(\widetilde{R}_{ij}\Phi,\omega,1,j(d-j+1))$$
$$+ \sum_{d_1+2\leq i\leq d} \Sigma_{\mathrm{st},i,1}(R_{d_1i+1}\Phi,\omega,1,(i+1)(d-i))$$
$$+ \sum_{d_1+2\leq i\leq d} \Sigma_{\mathrm{st},i,2}(\widetilde{R}_{d_1i}\Phi,\omega,1,i(d-i+1)).$$

By (4.5) and the propositions in this section, we get the following proposition.

**Proposition (5.9)**
(1) If $d$ is odd,

$$I(\Phi,\omega,w) \sim \delta(\omega_1)\delta(\omega_2)(\widehat{\Phi}(0) - \Phi(0))\Lambda(w;\rho)$$
$$+ \sum_{i=1,2}(F_i(\widehat{\Phi},\omega^{-1}) - F_i(\Phi,\omega))\varrho\Lambda(w;\rho).$$

(2) If $d = 2d_1$,

$$I(\Phi,\omega,w) \sim \delta(\omega_1)\delta(\omega_2)(\widehat{\Phi}(0) - \Phi(0))\Lambda(w;\rho)$$
$$+ \sum_{i=1,2,3}(F_i(\widehat{\Phi},\omega^{-1}) - F_i(\Phi,\omega))\varrho\Lambda(w;\rho)$$
$$+ \frac{\delta(\omega_1)\delta(\omega_2)}{(d_1+2)(d_1-1)}(R_{d_1d_1+1}\widehat{\Phi}(0) - R_{d_1d_1+1}\Phi(0))\varrho\Lambda(w;\rho)$$
$$+ \frac{1}{2}(T^1(\Phi,\omega) - T^1(\widehat{\Phi},\omega^{-1}))\varrho\Lambda(w;\rho).$$

Note that $\delta(\omega_2\omega_1^{-d_1})\delta(\omega_1) = \delta(\omega_1)\delta(\omega_2)$.

**Corollary (5.10)**
(1) If $d$ is odd,
$$I(\Phi,\omega) = \delta(\omega_1)\delta(\omega_2)\mathfrak{V}_2(\widehat{\Phi}(0) - \Phi(0))$$
$$+ \sum_{i=1,2}(F_i(\widehat{\Phi},\omega^{-1}) - F_i(\Phi,\omega)).$$

(2) If $d = 2d_1$,

$$I(\Phi, \omega) = \delta(\omega_1)\delta(\omega_2)\mathfrak{V}_2(\widehat{\Phi}(0) - \Phi(0))$$
$$+ \sum_{i=1,2,3} (F_i(\widehat{\Phi}, \omega^{-1}) - F_i(\Phi, \omega))$$
$$+ \frac{\delta(\omega_1)\delta(\omega_2)}{(d_1+2)(d_1-1)}(R_{d_1 d_1+1}\widehat{\Phi}(0) - R_{d_1 d_1+1}\Phi(0))$$
$$+ \frac{1}{2}(T^1(\Phi, \omega) - T^1(\widehat{\Phi}, \omega^{-1})).$$

## §6 The main statement

Let $[\frac{d}{2}] \leq i < j \leq d+1$. Then the following relations are easy to verify and the proof is left to the reader.

(6.1)
$$\widetilde{R}_{ij}\Phi_\lambda(x, x_j) = \lambda^{-(d-j+1)}\widetilde{R}_{ij}\Phi(\underline{\lambda}x, \underline{\lambda}^{-1}x_j),$$
$$\widetilde{R}_{ij}\widehat{\Phi}_\lambda(x, x_j) = \lambda^{-j}\widetilde{R}_{ij}\widehat{\Phi}(\underline{\lambda}^{-1}x, \underline{\lambda}x_j),$$
$$R_{ij}\Phi_\lambda(x) = \lambda^{-(d-j+1)}R_{ij}\Phi(\underline{\lambda}x),$$
$$R_{ij}\widehat{\Phi}_\lambda(x) = \lambda^{-j}R_{ij}\widehat{\Phi}(\underline{\lambda}^{-1}x).$$

The following proposition is an easy consequence of (6.1).

**Proposition (6.2)**
(1) Let $[\frac{d}{2}] < i$, $i+2 \leq j \leq d$. Then

$$\Sigma_{ij}(\widetilde{R}_{ij}\Phi_\lambda, \omega, 1, j(d-j+1)) = \lambda^{-(d-j+1)}\Sigma_{ij}(\widetilde{R}_{ij}\Phi, \omega, \lambda, j(d-j+1)),$$
$$\Sigma_{ij}(\widetilde{R}_{ij}\widehat{\Phi}_\lambda, \omega^{-1}, 1, j(d-j+1)) = \lambda^{-j}\Sigma_{ij}(\widetilde{R}_{ij}\widehat{\Phi}, \omega^{-1}, \lambda^{-1}, j(d-j+1)).$$

(2) Let $d = 2d_1$, $d_1 + 2 \leq i < j \leq d$. Then

$$\Sigma_{\text{st},ij}(\widetilde{R}_{d_1 j}\Phi_\lambda, \omega, 1, j(d-j+1)) = \lambda^{-(d-j+1)}\Sigma_{\text{st},ij}(\widetilde{R}_{d_1 j}\Phi, \omega, \lambda, j(d-j+1)),$$
$$\Sigma_{\text{st},ij}(\widetilde{R}_{d_1 j}\widehat{\Phi}_\lambda, \omega^{-1}, 1, j(d-j+1)) = \lambda^{-j}\Sigma_{\text{st},ij}(\widetilde{R}_{d_1 j}\widehat{\Phi}, \omega^{-1}, \lambda^{-1}, j(d-j+1)).$$

(3) Let $d = 2d_1$, $d_1 + 2 \leq i \leq d$. Then

$$\Sigma_{\text{st},i,1}(R_{d_1 i+1}\Phi_\lambda, \omega, 1, (i+1)(d-i))$$
$$= \lambda^{-(d-i)}\Sigma_{\text{st},i,1}(R_{d_1, j+1}\Phi, \omega, \lambda, (i+1)(d-i)),$$
$$\Sigma_{\text{st},i,1}(R_{d_1 i+1}\widehat{\Phi}_\lambda, \omega, 1, (i+1)(d-i))$$
$$= \lambda^{-(i+1)}\Sigma_{\text{st},i,1}(R_{d_1 i+1}\widehat{\Phi}, \omega, \lambda^{-1}, (i+1)(d-i)),$$
$$\Sigma_{\text{st},i,2}(\widetilde{R}_{d_1 i}\Phi_\lambda, \omega, 1, j(d-i-1))$$
$$= \lambda^{-(d-i+1)}\Sigma_{\text{st},i,2}(\widetilde{R}_{d_1 i}\Phi, \omega, \lambda, j(d-i-1)),$$
$$\Sigma_{\text{st},i,2}(\widetilde{R}_{d_1 i}\widehat{\Phi}_\lambda, \omega, 1, j(d-i-1))$$
$$= \lambda^{-i}\Sigma_{\text{st},i,2}(\widetilde{R}_{d_1 i}\widehat{\Phi}, \omega, \lambda^{-1}, j(d-i-1)).$$

**Proposition (6.3)** *The integrals*

$$\int_0^1 \lambda^s F_2(\widehat{\Phi}, \omega^{-1}) d^\times \lambda, \quad \int_0^1 \lambda^s F_3(\widehat{\Phi}, \omega^{-1}) d^\times \lambda$$

*converge absolutely and locally uniformly for all $s$.*

*Proof.* If $\lambda \leq 1$, $\lambda^{-1} \geq 1$. So by (2.2), (2.6), (2.9), and (2.12), for any $N \geq 1$, $\Sigma_{ij}(\widehat{\Phi_\lambda}, \omega^{-1}, 1, j(d-j+1))$ etc. are bounded by a constant multiple of of $\lambda^N$. This proves Proposition (6.3). □

**Proposition (6.4)** *Let $[\frac{d}{2}] < i$, $i + 2 \leq j \leq d$. Then the integral*

$$\int_0^1 \lambda^s \Sigma_{ij}(\widetilde{R}_{ij}\Phi_\lambda, \omega, 1, j(d-j+1)) d^\times \lambda$$

*converges absolutely and locally uniformly for $\operatorname{Re}(s) > d - i + 1$ if $d$ is odd, $j = d$, and $i = \frac{d+1}{2}$, and for $\operatorname{Re}(s) > d - i$ otherwise.*

*Proof.* By (2.2) and (6.2), for any $N \geq 1$,

$$\Sigma_{ij}(\widetilde{R}_{ij}\Phi_\lambda, \omega, 1, j(d-j+1)) \ll \lambda^{-(d-i-1+N)}$$

as long as

$$(2i - d)N - (2j - d) + j(d - j + 1) - 2 > 0.$$

It is easy to see that

$$-(2j - d) + j(d - j + 1) - 2 \begin{cases} = -2 & j = d, \\ > -1 & j < d. \end{cases}$$

Suppose $j < d$. Since $2i - d \geq 1$, $(2i - d) - (2j - d) + j(d - j + 1) - 2 > 0$. So we choose $N = 1$ in this case. Suppose $j = d$. If $2i - d \geq 2$, we choose $N = 1 + \delta$ where $\delta > 0$ is a constant. If $d$ is odd and $i = \frac{d+1}{2}$, we choose $N = 2 + \delta$ where $\delta > 0$ is a constant. Then $(2i - d)N - (2j - d) + j(d - j + 1) - 2 > 0$. So

$$\Sigma_{ij}(\widetilde{R}_{ij}\Phi_\lambda, \omega, 1, j(d-j+1)) \ll \begin{cases} \lambda^{-(d-i+1+\delta)} & d \text{ is odd, } j = d, \text{ and } i = \frac{d+1}{2}, \\ \lambda^{-(d-i+\delta)} & \text{otherwise.} \end{cases}$$

Since this inequality holds for all $\delta > 0$, we get Proposition (6.4). □

**Proposition (6.5)** *Let $d = 2d_1$, $d_1 + 2 \leq i < j \leq d$. Then the integral*

$$\int_0^1 \lambda^s \Sigma_{\operatorname{st},ij}(\widetilde{R}_{d_1 j}\Phi_\lambda, \omega, 1, j(d-j+1)) d^\times \lambda$$

*converges absolutely and locally uniformly for $\operatorname{Re}(s) > d - i + 3$.*

*Proof.* We choose $N = 1$ in (2.6). Then since $(2i-d)-(2j-d)+j(d-j+1)+2-2 > 0$ for all $d_1 + 2 \leq i < j \leq d$,

$$\Sigma_{\text{st},ij}(\widetilde{R}_{d_1 j}\Phi_\lambda, \omega, 1, j(d-j+1)) \ll \lambda^{-(d-i+3)}.$$

This proves Proposition (6.5). □

By a similar argument and using (2.9), we get the following proposition.

**Proposition (6.6)** *If $d = 2d_1$, $d_1 + 2 \leq i \leq d$, the integral*

$$\int_0^1 \lambda^s \Sigma_{\text{st},i,1}(R_{d_1 i+1}\Phi_\lambda, \omega, 1, (i+1)(d-i)) d^\times \lambda$$

*converges absolutely and locally uniformly for $\operatorname{Re}(s) > d - i + 4$.*

By (2.12), we get the following proposition.

**Proposition (6.7)** *If $d = 2d_1$, $d_1 + 2 \leq i \leq d$, the integral*

$$\int_0^1 \lambda^s \Sigma_{\text{st},i,2}(\widetilde{R}_{d_1 i}\Phi_\lambda, \omega, 1, i(d-i-1)) d^\times \lambda$$

*converges absolutely and locally uniformly for $s$.*

The following proposition follows from Propositions (6.4)–(6.7).

**Proposition (6.8)**
(1) *If $d$ is odd, the integral*

$$\int_0^1 \lambda^s F_2(\Phi_\lambda, \omega) d^\times \lambda$$

*converges absolutely and locally uniformly for $\operatorname{Re}(s) > \frac{d+1}{2}$.*
(2) *If $d$ is even, the integrals*

$$\int_0^1 \lambda^s F_2(\Phi_\lambda, \omega) d^\times \lambda, \quad \int_0^1 \lambda^s F_3(\Phi_\lambda, \omega) d^\times \lambda$$

*converge absolutely and locally uniformly for $\operatorname{Re}(s) > \frac{d}{2} + 2$.*

Let

$$p_i = d - i + \frac{(i+1)(d-i) - 2}{d - 2i}$$

for $[\frac{d}{2}] < i \leq d$. Then by (6.2),

$$\Sigma_{S_i}(\Phi_\lambda, \omega) = \lambda^{-p_i} \Sigma_{S_i}(\Phi, \omega),$$
$$\Sigma_{S_i}(\widehat{\Phi}_\lambda, \omega) = \lambda^{-(d+1-p_i)} \Sigma_{S_i}(\widehat{\Phi}, \omega).$$

Therefore, we get the following proposition.

**Proposition (6.9)**

(1) $$\int_0^1 \lambda^s \Sigma_{S_i}(\Phi_\lambda, \omega) d^\times \lambda = \frac{\Sigma_{S_i}(\Phi, \omega)}{s - p_i}.$$

(2) $$\int_0^1 \lambda^s \Sigma_{S_i}(\widehat{\Phi_\lambda}, \omega^{-1}) d^\times \lambda = \frac{\Sigma_{S_i}(\widehat{\Phi}, \omega^{-1})}{s - (d + 1 - p_i)}.$$

Now we consider the last two terms in (5.10).
By (6.1),

(6.10) $$\int_0^1 \lambda^s (R_{d_1 d_1 + 1}\widehat{\Phi_\lambda}(0) - R_{d_1 d_1 + 1}\Phi_\lambda(0)) d^\times \lambda$$
$$= \frac{R_{d_1 d_1 + 1}\widehat{\Phi}(0)}{s - 1 - d_1} - \frac{R_{d_1 d_1 + 1}\Phi(0)}{s - d_1}.$$

It is easy to see that

(6.11) $$\int_0^1 \lambda^s (T^1(\Phi_\lambda, \omega) - T^1(\widehat{\Phi_\lambda}, \omega^{-1})) d^\times \lambda$$
$$= -T_+(\Phi, \omega, s) - T_+(\widehat{\Phi}, \omega^{-1}, d + 1 - s) + T(\Phi, \omega, s).$$

By (3.13), $T(\Phi, \omega, s)$ is holomorphic for $\text{Re}(s) > 2$.
If $d \geq 4$, $\frac{d}{2} + 2 > 2$. So by (6.8), (6.10), and (6.11), we get the following theorem.

**Theorem (6.12)**
(1) If $d$ is odd,

$$Z(\Phi,\omega,s) = Z_+(\Phi,\omega,s) + Z_+(\widehat{\Phi},\omega^{-1},d+1-s)$$
$$+ \delta(\omega_1)\delta(\omega_2)\mathfrak{V}_2\left(\frac{\widehat{\Phi}(0)}{s-(d+1)} - \frac{\Phi(0)}{s}\right)$$
$$+ \sum_{[\frac{d}{2}]<i\leq d}\left(\frac{\Sigma_{S_i}(\widehat{\Phi},\omega^{-1})}{s-(d+1-p_i)} - \frac{\Sigma_{S_i}(\Phi,\omega)}{s-p_i}\right)$$
$$+ \int_0^1 \lambda^s(F_2(\widehat{\Phi_\lambda},\omega^{-1}) - F_2(\Phi_\lambda,\omega))d^\times\lambda,$$

and the last term is holomorphic for $\operatorname{Re}(s) > \frac{d+1}{2}$.
(2) If $d$ is even,

$$Z(\Phi,\omega,s) = Z_+(\Phi,\omega,s) + Z_+(\widehat{\Phi},\omega^{-1},d+1-s)$$
$$+ \delta(\omega_1)\delta(\omega_2)\mathfrak{V}_2\left(\frac{\widehat{\Phi}(0)}{s-(d+1)} - \frac{\Phi(0)}{s}\right)$$
$$+ \sum_{[\frac{d}{2}]<i\leq d}\left(\frac{\Sigma_{S_i}(\widehat{\Phi},\omega^{-1})}{s-(d+1-p_i)} - \frac{\Sigma_{S_i}(\Phi,\omega)}{s-p_i}\right)$$
$$+ \sum_{i=2,3}\int_0^1 \lambda^s(F_i(\widehat{\Phi_\lambda},\omega^{-1}) - F_i(\Phi_\lambda,\omega))d^\times\lambda$$
$$+ \frac{\delta(\omega_1)\delta(\omega_2)}{(d_1+2)(d_1-1)}\left(\frac{R_{d_1d_1+1}\widehat{\Phi}(0)}{s-1-d_1} - \frac{R_{d_1d_1+1}\widehat{\Phi}(0)}{s-d_1}\right)$$
$$- \frac{1}{2}\left(T_+(\Phi,\omega,s) + T_+(\widehat{\Phi},\omega^{-1},d+1-s) - T(\Phi,\omega,s)\right).$$

and the last three terms are holomorphic for $\operatorname{Re}(s) > \frac{d}{2} + 2$.

We list a few examples of the location of the poles we got as follows.

(1) If $d = 4$, $Z(\Phi,\omega,s)$ is meromorphic for $\operatorname{Re}(s) > 4$, and the poles in this domain are $s = \frac{9}{2}, 5$.
(2) If $d = 5$, $Z(\Phi,\omega,s)$ is meromorphic for $\operatorname{Re}(s) > 3$, and the poles in this domain are $s = 6, \frac{28}{5}, 10$.
(3) If $d = 6$, $Z(\Phi,\omega,s)$ is meromorphic for $\operatorname{Re}(s) > 5$, and the poles in this domain are $s = \frac{20}{3}, 7, 9$.
(4) If $d = 7$, $Z(\Phi,\omega,s)$ is meromorphic for $\operatorname{Re}(s) > 4$, and the poles in this domain are $s = \frac{54}{7}, 8, \frac{28}{3}, 18$.

Akihiko Yukie
Oklahoma State University
Mathematics Department
401 Math Science
Stillwater OK 74078-1058 USA
yukie@math.okstate.edu
http://www.math.okstate.edu/~yukie